%% file: HASEL_parallel_robot.tex
\documentclass{ifacconf}

\usepackage{graphicx}      
\usepackage{natbib}        

\usepackage{amsmath}
\usepackage{amssymb}
\usepackage{tikz}
\begin{document}
\begin{frontmatter}

\title{Design and Modeling of a HASEL Actuator-Based Micro Parallel Robot \thanksref{footnoteinfo}} 

\thanks[footnoteinfo]{This work is supported by the EIPHI Graduate School (contract: ANR-17-EURE-0002) and Région Franche-Comté.}

\author[First]{Agustin Feregrino}
\author[First]{Nelson Cisneros}
\author[First]{Alexis Lefèvre}
\author[First]{Yongxin Wu} 
\author[First]{Yann Le Gorrec} 

\address[First]{Université Marie et Louis Pasteur, SUPMICROTECH, CNRS, institut FEMTO-ST, F- 25000 Besançon, France (e-mails: yongxin.wu@femto-st.fr).}

\begin{abstract}
This paper presents the mechatronic design, dynamic modeling, and experimental validation of a three-degree-of-freedom (3-DOF) micro parallel robot featuring a prismatic–spherical (3PS) topology actuated by three Hydraulically Amplified Self-Healing Electrostatic (HASEL) actuators. Each soft actuator provides the prismatic motion of an individual limb, while a compliant interface to the moving platform functions as a spherical joint. A prototype incorporating three base-integrated HASEL actuators was fabricated, and the platform motion was measured using an XY laser-tracking system. For control purposes, a port-Hamiltonian (PH) model, combined with the mechanism’s forward kinematics (FKM), is developed to capture the robot’s nonlinear dynamic behavior, whereas the inverse kinematics (IKM) is employed to estimate the required actuator displacements. Model parameters were identified using nonlinear grey-box (NLGB) estimation, yielding a compact and control-oriented representation suitable for subsequent controller design.
\end{abstract}


\end{frontmatter}

\section{Introduction}
Soft robotics offers inherently safe, adaptable, and light weight systems, outperforming rigid robots in human interaction and unstructured environments \citep{KiLaTri:2013}. These capabilities are primarily enabled by the use of soft actuators, which provide compliant motion through flexible materials \citep{https://doi.org/10.1002/aisy.202000128}. Among soft actuators, Hydraulically Amplified Self-Healing Electrostatic (HASEL) actuators stand out for their large, controllable deformations, low energy consumption, and ease of fabrication, making them suitable for compact and portable robots \citep{MiWaAcMaLyKeGo:2019, Acome:2018}.
Their integration into different soft robotic devices is rapidly advancing due to their softness, precision, and efficiency. Recent studies highlight their use in biomedical diagnostics \citep{Ehrlich2025}, artificial sphincters \citep{lyko2024vitro}, and bio-inspired soft robots \citep{Perera2024, Xiong2024}. These developments highlight the need for accurate models capturing the complex dynamics of HASEL actuators.

Several modeling approaches for HASEL actuators have been reported in the literature, 
ranging from reduced-order lumped models \citep{Hainsworth2022} to data-driven approximations such as Dynamic Mode Decomposition \citep{Volchko2022}. While these methods are useful to approximate certain aspects of actuator behavior, they often neglect key nonlinear phenomena (such as drift effects and electromechanical coupling) or simplify the underlying electrical dynamics. As a result, their predictive capability is limited, particularly for control-oriented applications where an accurate representation of multiphysics interactions is essential \citep{Cisneros2024}.

The port-Hamiltonian (PH) framework offers a structured and physically interpretable approach to modeling systems composed of coupled physical domains by explicitly characterizing how energy is stored, dissipated, and exchanged among subsystems. This energy-based formulation makes PH modeling particularly well suited for describing HASEL actuators, as recently demonstrated in \citep{Cisneros2025}. While existing studies have primarily examined the dynamics of single HASEL actuators, the extension of PH modeling to mechanisms driven by multiple soft actuators remains largely unexplored. Such an extension is crucial for enabling coordinated motion generation, distributed actuation strategies, and systematic control synthesis in soft robotic architectures. Moreover, the intrinsic passivity properties of PH systems provide a strong foundation for control design, facilitating stable interconnections, and energy-aware control schemes.

Beyond its formal modeling advantages, the motivation for this work is also driven by the growing interest in compliant systems whose mechanical properties contribute directly to functional performance. This is particularly relevant at micro scales, where applications such as minimally invasive medical procedures, micro-manipulation, and the handling of delicate materials impose stringent requirements on precision, adaptability, and safe interaction. Fully exploiting these advantages requires control-oriented models capable of capturing the nonlinear interactions that arise in assemblies driven by multiple soft actuators. In this context, port-Hamiltonian (PH) formulations offer a natural framework, as they provide energy-consistent representations and preserve passivity, which is essential for achieving stable and predictable closed-loop behavior.

To investigate these considerations, this work examines the design and PH-based modeling of a three-degree-of-freedom (3-DOF) system actuated by three independent Hydraulically Amplified Self-Healing Electrostatic (HASEL) devices. Each actuation unit produces prismatic motion and is coupled to the moving structure through a passive spherical interface. This arrangement creates a fully HASEL-driven prismatic–spherical (3PS) configuration in which the nonlinear actuation characteristics interact with the geometric constraints of the mechanism.

The main contributions of this paper are:
\begin{itemize}
    \item We designed and fabricated a 3-DOF parallel platform actuated by HASEL devices, using additive manufacturing to obtain a lightweight and compliant prototype.
    \item We integrated the PH model of individual HASEL actuators, as reported in \cite{Cisneros2025}, into the kinematic framework of the parallel platform, enabling a unified description of the system.
    \item We validated the combined model experimentally by actuating the platform and tracking the end-effector motion.
    \item We provide a foundational framework for the design and analysis of soft parallel robots, validated through the first fully HASEL-actuated prototype. Our unified PHS model is a crucial step toward developing robust control strategies for this new class of compliant systems.
\end{itemize}

This paper is organized as follows. Section \ref{sec:design} presents the design and fabrication of the HASEL actuators and the 3PS parallel platform. Section \ref{sec:Model}  introduces the kinematic modeling of the mechanism and its connection with the PH actuator model. Section \ref{sec:identification} reports the parameter identification and experimental validation. Finally, Section \ref{sec:conclusion} concludes the work and outlines future directions.


\section{Design and manufacturing}\label{sec:design}
In this section, we describe the design and fabrication of the HASEL actuators and their integration into a 3-DOF parallel platform. 
We first detail the actuator manufacturing process, including the choice of materials and sealing methods.
We then present the design of the parallel platform.

\subsection{Micro Parallel Robot Design and Assembly}

The parallel platform is formed by a triangular structure supported by three connecting links, each attached at its base to a HASEL actuator. The kinematic behavior corresponds to a 3PS configuration. In this topology, the prismatic joints (P) are provided by the vertical deformation of each HASEL actuator, while the spherical joints (S) arise at the lower connection between each link and its respective actuator. This spherical behavior does not result from a mechanical ball joint but rather from the intrinsic flexibility of the links combined with the compliant attachment of the actuator to the rigid component.

\begin{figure}[ht!]
\centering
\resizebox{0.5\textwidth}{!}{\input{robotHASELsch_v2}} 
\caption{a) Top view HASEL actuator schematic design and 3PS parallel platform. b) Kinematic diagram of the 3PS topology. c) Side views showing the actuator displacements $h_i$, centroid $C_p$, and antenna $Tip$ with length $L$. d) HASEL schematic.} 
\label{fig:HASEL-schematic}
\end{figure}
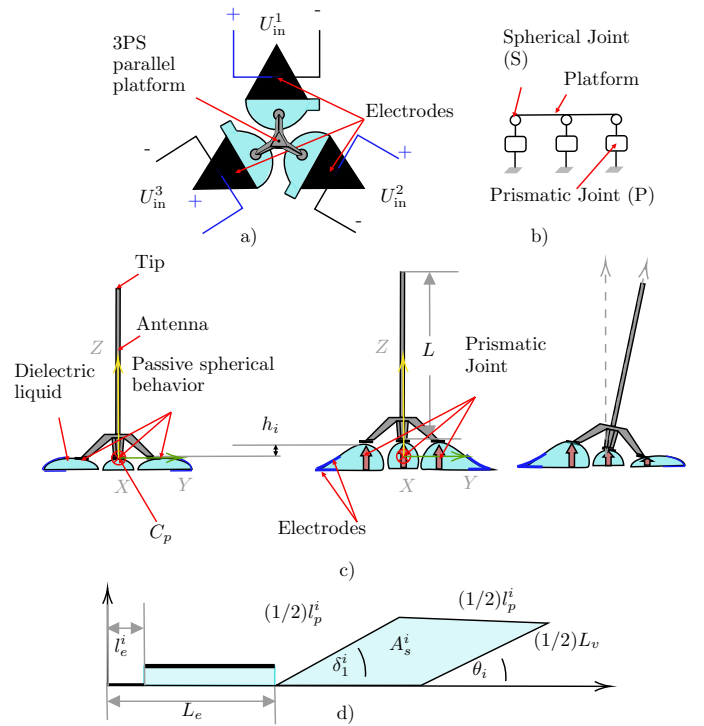

Figure~\ref{fig:HASEL-schematic} a) illustrates the design concept of the HASEL-actuated parallel platform. Figure~\ref{fig:HASEL-schematic} (b) shows the kinematic diagram of the 3PS topology, where each HASEL actuator provides the P and the flexible link connections emulate S. Figure~\ref{fig:HASEL-schematic} (c) depicts the schematic representation of the 3D-printed prototype in its neutral position, highlighting the antenna for motion tracking, the dielectric shell, and the compliant electrodes and the actuated configuration: when a high voltage $U_{in}$ is applied, the resulting prismatic deformation of the actuators (red arrows) induces platform motion through passive spherical behavior.


\subsection{Fabrication of the Actuators}

The actuators were fabricated using thermoplastic BOPP film as the dielectric shell, filled with Envirotemp FR3 biodegradable transformer oil as the dielectric fluid. Carbon-based conductive paint was applied as flexible electrodes on the actuator surfaces \citep{Acome:2018, Rothemund2020}. The manufacturing protocol was adapted from the procedure reported in \citep{Cis:2025}, including cutting, heat sealing, filling, and assembly.

The detailed manufacturing procedure consists of the following steps:

\begin{enumerate}
    \item \textbf{Sealing of the outline.}  
    The actuator outlines were sealed using a 3D printer operated without filament. The heat from the printer nozzle was applied directly to the BOPP film to create a hermetic perimeter.
    
    \item \textbf{Filling with dielectric fluid.}  
    Each actuator was filled with 0.4\,mL of dielectric fluid using a precision syringe. The fluid was manually redistributed to minimize air bubbles.
    
    \item \textbf{Application of compliant electrodes.}  
    Compliant electrodes were applied using a template and allowed to dry before final sealing.
    
    \item \textbf{Final sealing.}  
    The perimeter was sealed with a soldering iron at 217\,$^\circ$C, which provided reliable sealing without puncturing or dielectric breakdown.
    
\end{enumerate}

\subsection{Fabrication of the platform}

The structural components were fabricated using Fused Deposition Modeling (FDM) with PLA filament. 
Each print required approximately 23 minutes, enabling rapid prototyping and straightforward replacement during experimental trials.

\begin{figure}[ht!]
\centering
\includegraphics[width=4cm]{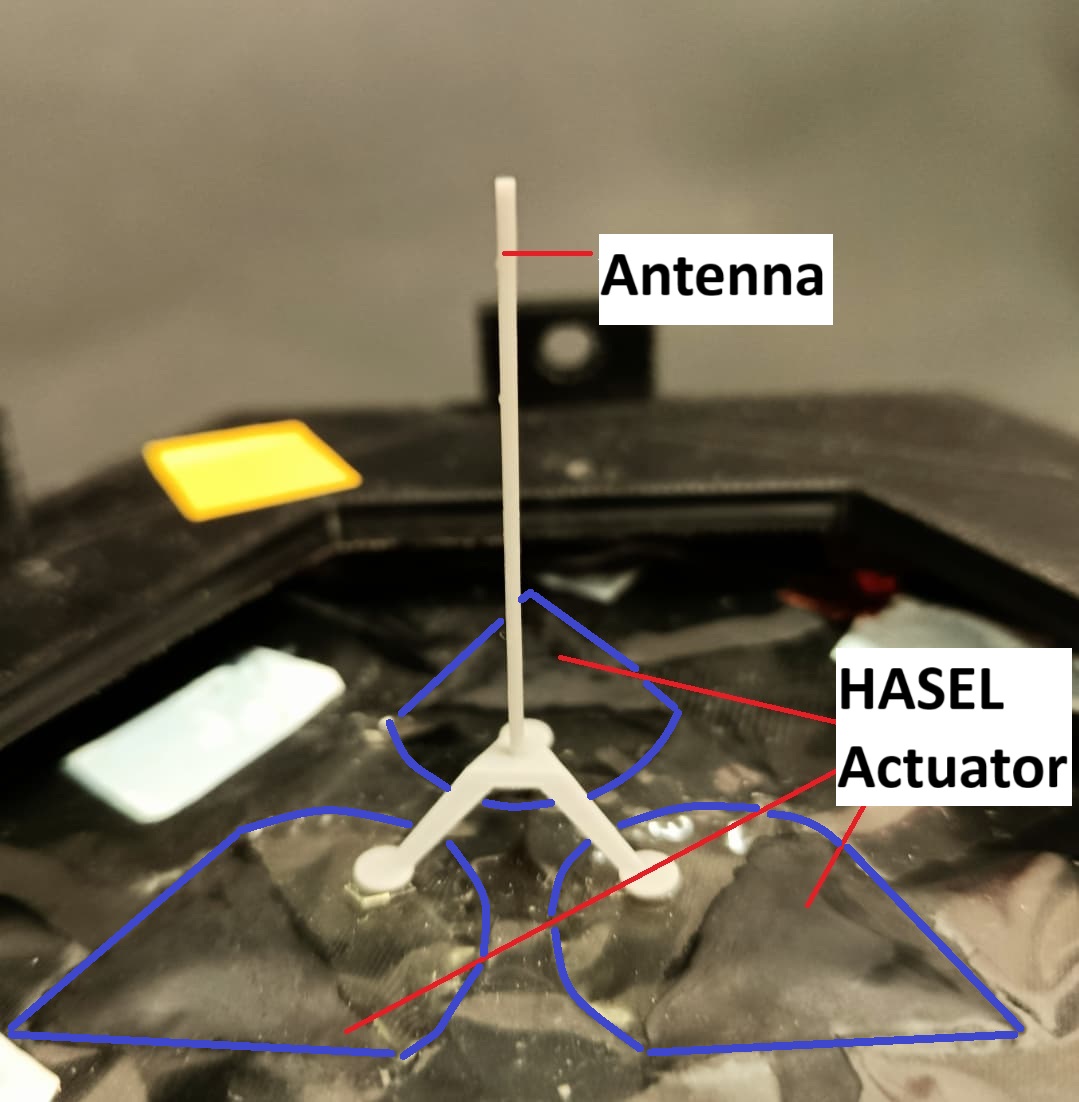}
\caption{ 3PS prototype showing the integrated antenna for motion tracking and the three HASEL actuators providing the prismatic inputs.}
\label{fig:prototype-platform}
\end{figure}

Figure~\ref{fig:prototype-platform} shows the physical prototype of the platform after integration of the HASEL actuators. The actuators provide the prismatic inputs to the 3PS mechanism.

\section{Modeling of micro parallel robot}\label{sec:Model}

In this section we present the soft parallel robot model driven by HASEL actuators based on the HASEL actuator PH model from \citep{Cisneros2025} and  the kinematic analysis of the 3PS parallel platform.  The model establishes the relation between actuator elongations and the motion of the antenna. 
Forward Kinematics Model (FKM) computes the antenna tip position from the actuator displacements, while Inverse Kinematics Model (IKM) determines the actuator heights for a desired tip position. 

\subsection{Control oriented port-Hamiltonian formulation}
The PH formalism provides a general framework to model multi-domain dynamical systems in terms of energy storage, dissipation, and interconnection through power ports \citep{Ortega2001,vdSchaft2002}. It can be expressed as:
\vspace{-0.2cm}
\begin{equation}
\begin{array}{ll}
\\& \dot{x} = [J(x) - R(x)] \nabla H(x) + G(x)u, 
\\& y = G^\top(x)\nabla H(x),
\end{array}
\end{equation}
where $x \in \mathbb{R}^n$ is the state vector, $H(x)$ the Hamiltonian function representing stored energy, $J(x)=-J^\top(x)$ the interconnection matrix, $R(x)=R^\top(x)\geq 0$ the dissipation, and $G(x)$ the input mapping. The output $y$ corresponds to the conjugate power variables, ensuring energetic consistency.


HASEL actuators combine electrostatic forces, soft dielectric shells, and internal fluids to produce large deformations with high efficiency \citep{Rothemund2020}. 

The total energy of the system is expressed as: 
\begin{equation}
    \begin{array}{ll}
    H=&\frac{1}{2}\theta^\top K_b\theta+ H_{\bar{g}} +\frac{1}{4}(l_{p}-L_p)^\top K(l_{p}-L_p)\\&+ \frac{1}{2}p^\top M^{-1} p  +\frac{1}{2}Q_1^\top C_1^{-1}Q_1+\frac{1}{2}Q_2^\top C_2^{-1}(\theta,l_p)Q_2
    \end{array}
\end{equation}

The angular deformation is $\theta=[\theta_1 \; \theta_2 \;\theta_3]^{\top}$, the top film deformation is $l_p=[l^1_p \; l^2_p \; l^3_p]^\top$, the angular momentum is $p=[p_1 \; p_2 \; p_3]^\top$. The electrical charges are defined as $Q_1=[Q^1_1 \; Q^2_1 \;Q^3_1]^\top$ and $Q_2=[Q^1_2 \; Q^2_2 \;Q^3_2]^\top$.
$M$ is the inertia matrix defined as $M=\text{diag}([\Tilde{I_1} \; \Tilde{I_2} \; \Tilde{I_3}])$ where the inertia is $\Tilde{I_i}=(1/24)L_v m$. The mass of each HASEL actuator is $m$. The potential energy related with the gravity is defined as $H_{\bar{g}}=1/2\bar{g}L_vm(\sin{\theta_1}+sin{\theta_2}+sin{\theta_3})$. $L_v$ is the bottom film length.

 $C_1=\text{diag} ([C^1_{1} \;C^2_{1} \; C^3_{1}])$ is the constant capacitance that represents the actuator's charge retention effect. $C_2=\text{diag} [C^1_{2} \;C^2_{2} \; C^3_{2}]$ is the dynamic capacitance that represents the electrodes, films, and dielectric liquid, computed as the sum of the capacitance of the zipped electrodes and the capacitance of the unzipped part of the electrodes.  
\begin{align}
    C^i_{2}(\theta_i,l_p^i)=\epsilon_0 \epsilon_r w\left( \frac{l^i_{e}(\theta_i,l_p^i)}{2t}+\frac{ L_e-l^i_{e}(\theta_i,l_p^i)}{2t+X_h}\right)
\end{align}
where $\epsilon_0$ is the vacuum permittivity, $\epsilon_r$ is the relative permittivity, $w$ is the actuator width, $t$ is the film thickness, $L_e$ is the electrode's length and $l_e^i$ is the length of the zipped electrodes part,  where $i=1,2,3$.

The zipped electrodes length is:

 \begin{equation}\label{eq:le}
      l^i_{e}(\theta_i,l_p^i)=L_e-\frac{1}{X_h}\left(A_T-\frac{L_vl^i_{p}}{4}sin (\delta^i_{1}(\theta_i,l_p^i))\right).
\end{equation}

Where $X_h$ is the distance between the unzipped electrodes. The constant area of the liquid inside the shell is $A_T$.

The dynamics of the system are:

\begin{equation}
\begin{array}{ll}
\underbrace{
\begin{bmatrix}

\dot{\theta} \\
\dot{l}_{p} \\
\dot{p} \\
\dot{Q}_{1} \\
\dot{Q}_{2}
\end{bmatrix}
}_{\dot{x}}
=&
\underbrace{
\begin{bmatrix}
0 & 0 & I & 0 & 0 \\
0 & 0 & d & 0 & 0 \\
-I & -d & -b & 0 & 0 \\
0 & 0 & 0 & -(\bar{R}_{0}+\bar{R}_{1}) & -\bar{R}_{0} \\
0 & 0 & 0 & -\bar{R}_{0} & -(\bar{R}_{0}+\bar{R}_{2})
\end{bmatrix}
}_{J-R}
\underbrace{
\begin{bmatrix}
\nabla_{\theta} H \\
\nabla_{l_{p}} H \\
\nabla_{p} H \\
\nabla_{Q_{1}} H \\
\nabla_{Q_{2}} H
\end{bmatrix}
}_{\nabla_{x} H} \\[1em]
&+
\underbrace{
\begin{bmatrix}  
0 \;0 \; 0 \; 
\bar{R}_{0} g a(\theta)^{\top} \;
\bar{R}_{0} g a(\theta)^{\top}
\end{bmatrix}^\top
}_{g} U_{in}
\label{eq:ph-system}
\end{array}
\end{equation}

\begin{equation}
y =
\underbrace{\big(\bar{R}_{0} g a(\theta)\big)^{\top} C_{1}^{-1} Q_{1}
+ \big(\bar{R}_{0} g a(\theta)\big)^{\top} C_{2}^{-1} Q_{2}}_{g^{\top}\nabla_{x}H},
\end{equation}
where $x=[{\theta}^\top,{l}_{p}^\top,{p}^\top,{Q}_{1}^\top,{Q}_{2}^\top]^{\top}$ includes all electrical and mechanical state variables. 

The coupling term $d$ is defined as $d=\text{diag}([\frac{2A^1_s}{l^1_p} \; \frac{2A^2_s}{l^2_p} \; \frac{2A^3_s}{l^3_p}])$. The area inside the shell is computed by:

\begin{equation}\label{eq_A2_0_2}
   A^i_{s}=\frac{1}{4}l^i_{p} L_v sin (\delta^i_{1}(\theta_i,l_p^i)). 
\end{equation}
 where 
 \begin{equation}\label{eq:delta_1}
      \delta^i_{1}(\theta_i,l_p^i)=\frac{\pi+\theta_i}{2}-\sin^{-1}\left(\frac{L_v}{l^i_{p}}\sin\left(\frac{\pi-\theta_i}{2}\right)\right)
\end{equation}
The damping of the system is $b=\text{diag}([b_1 \; b_2 \; b_3])$. The admittances of the system are $\bar{R}_{0}$, $\bar{R}_{1}$ and $\bar{R}_{2}$. They are defined as $\bar{R}_{0}=\text{diag}([\frac{1}{R^1_0} \; \frac{1}{R^2_0} \; \frac{1}{R^3_0}])$, $\bar{R}_{1}=\text{diag}([\frac{1}{R^1_1} \; \frac{1}{R^2_1} \; \frac{1}{R^3_1}])$ and $\bar{R}_{2}=\text{diag}([\frac{1}{R^1_2} \; \frac{1}{R^2_2} \; \frac{1}{R^3_2}])$. The term $ga$ is defined as $ga=\text{diag}([ga_1 \; ga_2\;ga_3])$, where $ga_i=\gamma_1 \cos{(\gamma_2 \theta_i)}$. $\gamma_1$ and $\gamma_2$ are parameters to be identified. 
The input voltage of the system is $U_{\text{in}}=[U^1_{\text{in}} \; U^2_{\text{in}} \; U^3_{\text{in}}]$. 

The vertical elongation of the HASEL actuator is:
\begin{equation}\label{eq:prismatic_motion}
    h_p=[\frac{1}{2}L_v \sin{\theta_1}, \;\; \;\; \frac{1}{2}L_v \sin{\theta_2},  \;\; \;\;\frac{1}{2}L_v \sin{\theta_3}]^\top
\end{equation}
This detailed formulation serves as the foundation for embedding the actuator into the kinematic framework of our 3-DOF platform.

\subsection{Forward Kinematics}\label{sec:fkm}

\begin{figure}[h]
\centering
\includegraphics[width=1\linewidth]{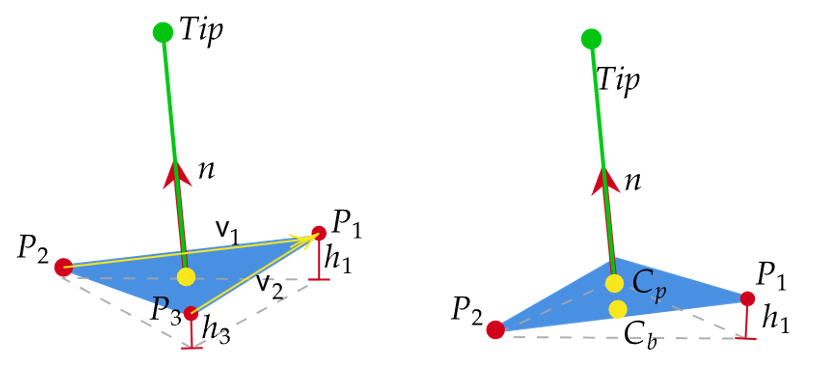}
\caption{Left: Forward kinematics; Right: Inverse kinematics}
\label{fig:fkm-ikm-geometry}
\end{figure}

To determine the spatial position of the end-effector ({\it i.e.,} the tip of the antenna) based on the vertical elongation of three HASEL actuators, the Forward Kinematics Method (FKM) is employed. 
Assume that each actuator is located at fixed planar coordinates $(X_i, Y_i)$, equally spaced at the vertices of an equilateral triangle. 
The vertical elongation of $i$-th actuator is denoted as $h_i(t)$. 
Accordingly, the position vector of the $i$-th platform corner point can be expressed as:
\begin{equation}
    \mathbf{P}_{i}(t)=\left[\begin{array}{c}
X_{i} \\
Y_{i} \\
h_{i}(t)
\end{array}\right].
\end{equation}
The centroid of the platform is given by:
\begin{equation}
    \mathbf{C_p}(t) = \frac{1}{3}\big(\mathbf{P}_1(t) + \mathbf{P}_2(t) + \mathbf{P}_3(t)\big)
\end{equation}
The orientation is obtained by defining the unit normal vector as:
\begin{equation} 
\begin{array}{l}
\mathbf{v}_1 = \mathbf{P}_2(t) - \mathbf{P}_1(t), \\ 
\mathbf{v}_2 = \mathbf{P}_3(t) - \mathbf{P}_1(t), \\   
\mathbf{n}(t) = \frac{\mathbf{v}_1 \times \mathbf{v}_2}{\|\mathbf{v}_1 \times \mathbf{v}_2\|} 
\end{array}
\end{equation}
Assuming the antenna is rigid and perpendicular to the platform, the tip position is:
\begin{equation}
    \mathbf{Tip}(t) = \mathbf{C_p}(t) + L \, \mathbf{n}(t)
    \label{eq:tip}
\end{equation}
where $L$ is the antenna length.
 The actuators are located at fixed base positions $P_i$, with vertical displacements $h_i(t)$ defining the orientation of the platform. 
The centroid $\mathbf{C_p}(t)$ and the unit normal vector $\mathbf{n}$ determine the position of the antenna tip $\mathbf{Tip}(t)$ according to \eqref{eq:tip} (left figure of Fig. \ref{fig:fkm-ikm-geometry}). The same representation will be used in the derivation of IKM (see right figure of Fig. \ref{fig:fkm-ikm-geometry}).

\subsection{Inverse Kinematics (normal-aligned stick)} \label{sec:ikm}

For the identification and control design purposes, it is necessary to compute the actuator heights $h_i$ at the three base locations $(X_i, Y_i)$ from the tip position $\mathbf{Tip} = \bigl[ X,\, Y,\, Z \bigr]^\top$, which is rigidly attached and remains orthogonal to the platform.

Let $\mathbf{C}_b=\big(\frac{1}{3}\sum_i X_i,\ \frac{1}{3}\sum_i Y_i,\ 0\big)$ denote the centroid of the base in the horizontal plane. 
A direction vector from $\mathbf{C}_b(t)$ to the tip position $\mathbf{Tip}(t)$ is defined as:
\begin{equation}
    \mathbf{v}=\left[\begin{array}{c}
X-\frac{1}{3} \sum_{i} X_{i} \\
Y-\frac{1}{3} \sum_{i} Y_{i} \\
Z
\end{array}\right], \quad \mathbf{n}=\frac{\mathbf{v}}{\|\mathbf{v}\|}=\left[\begin{array}{l}
n_{X} \\
n_{Y} \\
n_{Z}
\end{array}\right]
\end{equation}
which serves as an estimate of the platform normal, assumed to be aligned with the antenna (or stick) axis.
Because the stick is orthogonal to the platform and has length $L$, the platform centroid is located along $- \mathbf{n}$ from the tip:
\begin{equation}
    \mathbf{C_p}=\operatorname{\mathbf{Tip}}-L \mathbf{n}=\left[\begin{array}{l}
C_{pX} \;
C_{pY} \;
C_{pZ}
\end{array}\right]^\top
\end{equation}

The platform is modeled as a plane with normal $\mathbf{n}$ passing through $\mathbf{C_p}$, i.e.,
\begin{equation}
    n_X(X - C_{pX}) + n_Y(Y - C_{pY}) + n_Z(Z - C_{pZ}) = 0.
\end{equation}
Evaluating this relation at the base anchor $(X_i,Y_i)$ yields the platform height $h_i$ at each anchor:
\begin{equation}
    h_i \;=\; C_{pZ} \;-\; \frac{\,n_X (X_i - C_{pX}) + n_Y (Y_i - C_{pY})\,}{\,n_Z\,}
    \label{eq:IK-hi}
\end{equation}


To accurately capture the motion of the 3-DOF parallel platform, the kinematic model is enhanced by incorporating the physical dynamics of the HASEL actuators, as described by \eqref{eq:ph-system}.
In particular, we use the output corresponding to the vertical displacement \eqref{eq:prismatic_motion}, which characterizes the effective prismatic motion at joint $i$. 

By mapping the actuator displacement $h_i(t)$ to the associated platform corner, a direct relation is established between the control input (applied high voltage) and the geometric variables used in both the forward and inverse kinematics formulations.  This approach allows the consistent description of the platform motion from the electrical input to the resulting tip position of the antenna.

Furthermore, this integration highlights the modularity of the framework: the same PH actuator model can be embedded in multi-actuator systems without modification, thereby enabling predictive simulation and control of the complete 3-DOF platform. 
In addition, any nonlinear or energy-dependent effects captured by the PH formulation are naturally reflected in the platform kinematics, yielding a more physically consistent and realistic representation than that provided by purely geometric models.


\section{Identification and Validation}\label{sec:identification}

\subsection{Experimental Excitation and Data Acquisition}

\begin{figure}[h]
\centering
\includegraphics[width=1\linewidth]{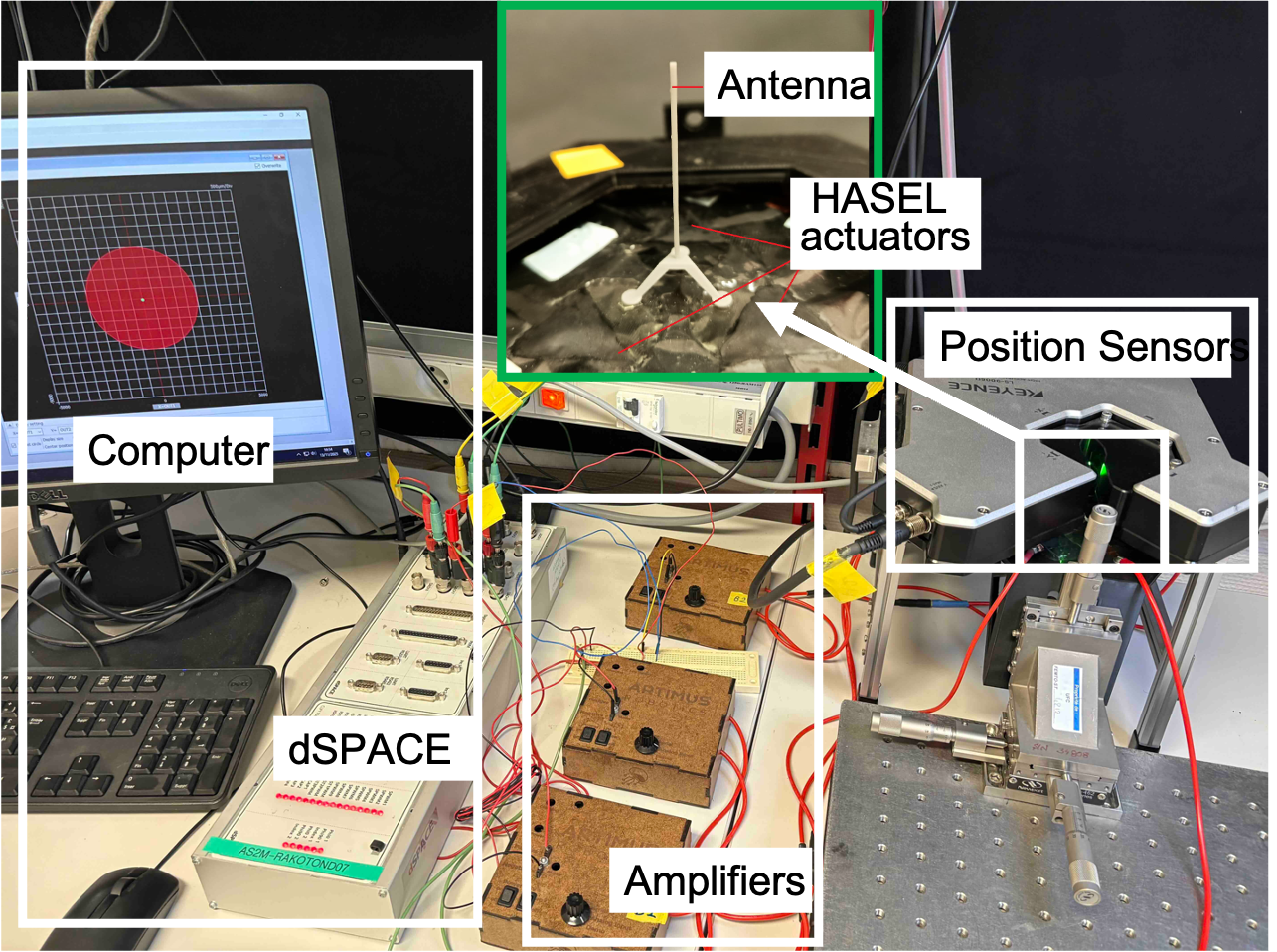}
\caption{Overall experimental setup}
\label{fig:experiemental seuup}
\end{figure}

To identify the parameters of each actuator, the three HASEL units were driven by sinusoidal high-voltage inputs of the form $
u_i(t) = U_0 \sin(\omega t + \varphi_i)$,
where the excitation frequency was set to $\omega = 3\pi \,\text{rad/s}$. The input signals were phase-shifted by $120^\circ$ relative to each other, generating an approximately circular trajectory of the platform tip.
The tip displacement was recorded using a high-speed 2D laser displacement sensor (Keyence LS-9000) with the sampling time of $0.1$ ms. 
The acquired measurements were subsequently processed using the inverse kinematic model (Section~\ref{sec:ikm}) to extract the effective vertical displacements $h_i(t)$ corresponding to each actuator. The input signal generation and the data acquisition are connected to the computer via a dSPACE card. The overall experimental setup is shown in Fig. \ref{fig:experiemental seuup}.

\subsection{Nonlinear Grey-Box Identification}

The dynamics of each HASEL actuator were represented using a nonlinear ordinary differential equation (ODE) model, as described by \eqref{eq:ph-system}.
The adopted grey-box model included nine parameters: (i) electrical parameters ($R_0, R_1, R_2$) representing equivalent resistances and $C_1$ denoting the effective capacitance;
(ii) mechanical parameters $K$,  $K_b$ and $b$ denoting the  linear stiffness, the torsional stiffness and the damping coefficient, respectively; 
(iii)  input gain coefficients $\gamma_1, \gamma_2$. 
This grouping preserves the physical interpretability of the model while allowing parameter estimation from experimental data.

Parameter identification was performed in MATLAB using the \texttt{System Identification Toolbox}. 
A nonlinear grey-box (NLGB) strategy was implemented via the \texttt{nlgreyest} routine, employing the Levenberg–Marquardt (LM) optimization algorithm. 
The solver \texttt{ode15s} was used to handle the stiffness of the actuator dynamics, with absolute and relative tolerances set to $10^{-3}$. 

\begin{figure}[h]
\centering
\includegraphics[width=1\linewidth]{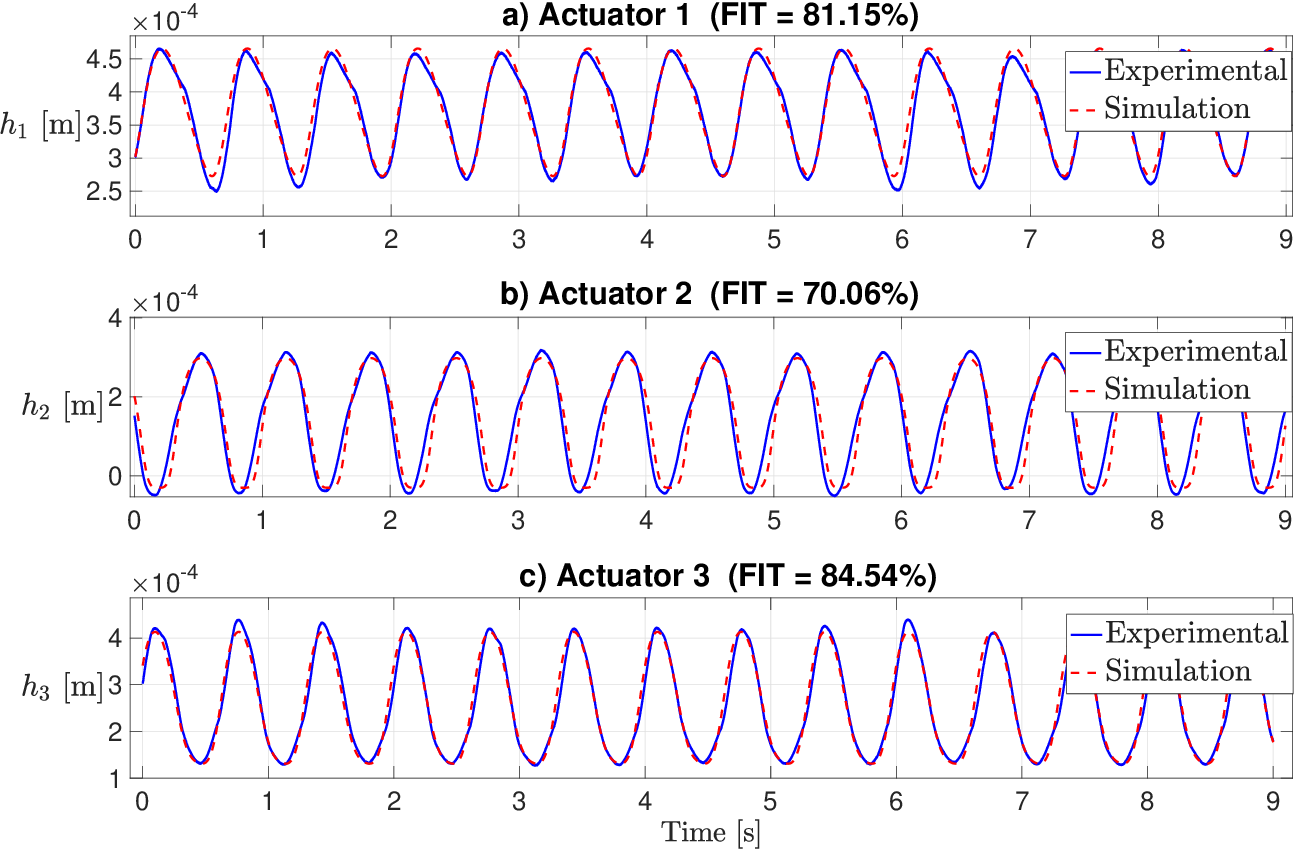}

\medskip

\includegraphics[width=0.5\linewidth]{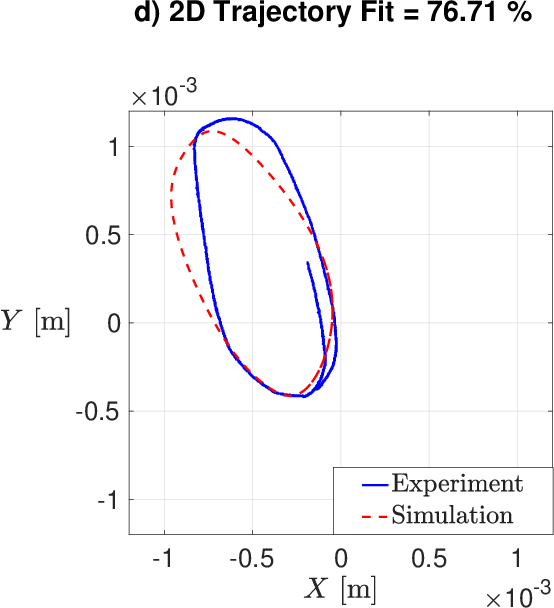}

\caption{Parameter identification results. 
(a–c) Experimental vs. simulated vertical displacements $h_i(t)$ for each actuator. 
(d) Reconstructed tip trajectory obtained from the identified models via FKM.}
\label{fig:identification-results}
\end{figure}

Figures~\ref{fig:identification-results} (a)–(c) compare the experimental vertical displacements $h_i(t)$ of each actuator ($i=1,2,3$) with the simulated response of the identified PH-based grey-box model. The  experimental displacement is found 
by measuring the position of the tip and using the IKM. The corresponding fit percentages, expressed as the normalized root-mean-square error (NRMSE), are reported in each subplot.
  
Figure~\ref{fig:identification-results} (d) shows the reconstructed tip trajectory in the XY plane, obtained by applying the forward kinematic model (Section~\ref{sec:fkm}) to the simulated actuator outputs; the associated 2D fit percentage is also indicated.  
These results confirm that the identified models accurately capture both the actuator-level dynamics and the coordinated motion of the complete platform.  
During identification at $\omega = 3\pi $\,rad/s, the models achieved fit percentages (NRMSE) of 81.15\%, 70.06\%, and 84.54\% for Actuators~1--3, respectively, demonstrating the high accuracy of the grey-box parameter identification.

The final identified values of the electrical, mechanical, and coupling parameters are summarized in Table~\ref{tab:identified_parameters}, which reports the effective capacitances, stiffness and damping coefficients, and input gains associated with each actuator.

\begin{table}[ht]
\centering
\caption{Model parameters of the three HASEL actuators.}
\label{tab:identified_parameters}
\begin{tabular}{lccc}
\hline
\textbf{Parameter} & \textbf{Act. 1} & \textbf{Act. 2} & \textbf{Act. 3} \\
\hline
$K$ [N/m]        & 0.250 & 0.200 & 0.250 \\
$K_b$ [Nm/rad]   & 0.305 & 0.300 & 0.300 \\
$b$ [kg$\cdot$s] & 0.129 & 0.010 & 0.015 \\
$R_{0}$ [$\Omega$] & 600 & 600 & 600 \\
$R_{1}$ [$\Omega$] & 99.98 & 112.7 & 100 \\
$R_{2}$ [$\Omega$] & 5.385$\times 10^{4}$ & 5.60$\times 10^{4}$ & 5.5$\times 10^{4}$ \\
$C_{1}$ [F]        & 2.15$\times 10^{-10}$ & 2.00$\times 10^{-10}$ & 2.20$\times 10^{-10}$ \\
$\gamma_{1}$ [–]   & 118.4 & 131.9 & 118.4 \\
$\gamma_{2}$ [–]   & 58.0 & 30.5 & 25.0 \\
\hline
\multicolumn{4}{c}{\textbf{Common parameters for all actuators}} \\
\hline
$L_{p}$ [m]       & \multicolumn{1}{c}{0.014} &
$L_{v}$ [m]       & \multicolumn{1}{c}{0.012} \\
$L_{e}$ [m]       & \multicolumn{1}{c}{0.020} &
$X_{h}$ [m]       & \multicolumn{1}{c}{0.002} \\
$m$ [kg]          & \multicolumn{1}{c}{0.001133} & 
$\epsilon_{r}$ [F/m] & \multicolumn{1}{c}{2.2} \\
$\epsilon_{0}$ [F/m] & \multicolumn{1}{c}    {$8.854\times 10^{-12}$} & 
$w$ [m]           & \multicolumn{1}{c}{0.025} \\
$t$ [m]           & \multicolumn{1}{c}{$18\times 10^{-6}$}        &
$g$ [m/s$^{2}$]   & \multicolumn{1}{c}{9.8} \\
\hline
\end{tabular}
\end{table}

\subsection{Model Validation}

To further evaluate the predictive capability of the identified models, a validation experiment was performed using an independent excitation frequency of $\omega = 11 \,\text{rad/s}$, distinct from the $3\pi \,\text{rad/s}$ dataset used for estimation. Each actuator was driven by sinusoidal voltages phase-shifted by $120^\circ$, and the resulting displacements were recorded with the 2D laser sensor.

\begin{figure}[h]
\centering
\includegraphics[width=1\linewidth]{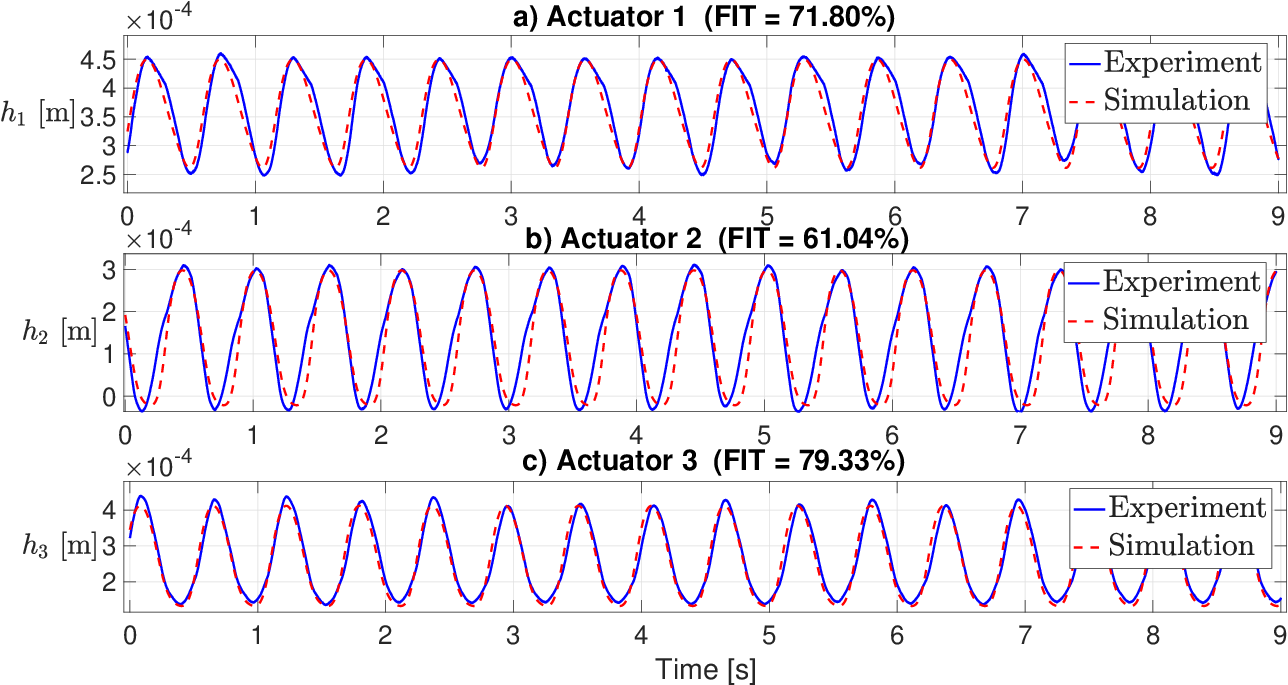}

\vspace{0.5em}

\includegraphics[width=0.5\linewidth]{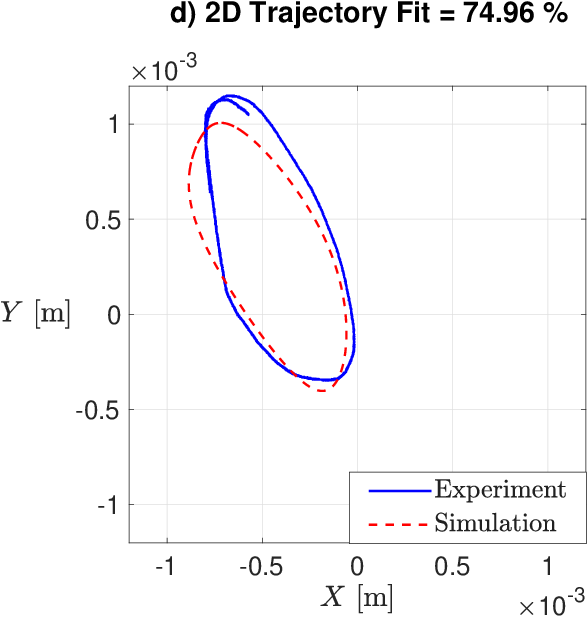}

\caption{Validation results under $\omega=11$\,rad/s. 
(a–c) Experimental vs. simulated vertical displacements $h_i(t)$; 
(d) Reconstructed tip trajectory obtained via FKM.}
\label{fig:validation-results}
\end{figure}

The simulated responses of the PH-based grey-box models showed good agreement in both amplitude and phase. The NRMSE fits achieved during validation were 71.80\%, 61.04\%, and 79.33\% for Actuators~1--3, respectively.  
The reconstructed tip trajectory in the XY plane, obtained via the forward kinematic model (Section~\ref{sec:fkm}), also followed the experimental motion with a 2D fit of 74.96\%, as shown in Fig.~\ref{fig:validation-results}.  
These results demonstrate that the identified parameters generalize effectively the training conditions, confirming that the nonlinear PH framework combined with grey-box identification captures the essential dynamics of the HASEL actuators and supports system-level modeling of the 3PS platform.

\section{Conclusion}\label{sec:conclusion}

In this paper we presented the design, modeling, and experimental validation of a three-degree-of-freedom soft parallel robot actuated by HASEL devices. The robot is based on a 3PS topology, where each limb combines one prismatic joint actuated by a HASEL actuator and one passive spherical joint provided by the flexible links.  

The PH framework was employed to describe the multiphysics dynamics of the actuators, and integrated into the kinematic formulation of the parallel mechanism. A nonlinear grey-box identification approach was used to estimate the physical parameters of the model. The identified PH models reproduced the actuator dynamics with high accuracy, achieving NRMSE fits above 70\% during both identification and validation. The reconstructed platform trajectory also matched experimental measurements, confirming the predictive capability of the proposed model.  

The perspectives of this work are to implement closed-loop control strategies based on the PH formulation, extend the modeling to capture multi-actuator coupling effects, and apply the framework to biomedical and microscale robotic systems.

\small

\bibliography{ifacpaper_2026}             
                                                   







\end{document}

%% file: robotHASELsch_v2.tex
\tikzset{every picture/.style={line width=0.75pt}} 

\begin{tikzpicture}[x=0.75pt,y=0.75pt,yscale=-1,xscale=1]

\draw  [draw opacity=0][fill={rgb, 255:red, 74; green, 221; blue, 226 }  ,fill opacity=0.5 ] (319.19,151.22) .. controls (310.3,141.92) and (307,130.36) .. (311.86,121.62) .. controls (316.73,112.85) and (328.46,109.32) .. (341.28,111.66) -- cycle ;
\draw  [draw opacity=0][fill={rgb, 255:red, 74; green, 221; blue, 226 }  ,fill opacity=0.5 ] (268.1,111.98) .. controls (281.5,108.87) and (293.84,111.85) .. (298.92,120.54) .. controls (303.86,128.99) and (300.57,140.56) .. (291.63,150.19) -- cycle ;
\draw  [draw opacity=0][fill={rgb, 255:red, 74; green, 221; blue, 226 }  ,fill opacity=0.5 ] (327.37,88.4) .. controls (323.4,100.97) and (314.73,109.69) .. (304.69,109.69) .. controls (294.61,109.69) and (285.93,100.92) .. (281.97,88.29) -- cycle ;
\draw    (325.22,95.69) .. controls (308.48,121.6) and (281.26,104.33) .. (283.36,87.88) ;
\draw  [fill={rgb, 255:red, 0; green, 0; blue, 0 }  ,fill opacity=1 ] (304.29,50.47) -- (325.22,87.88) -- (283.36,87.88) -- cycle ;
\draw    (325.22,87.88) -- (333.6,87.88) ;
\draw    (325.22,95.69) -- (333.6,95.69) ;
\draw  [draw opacity=0][fill={rgb, 255:red, 74; green, 221; blue, 226 }  ,fill opacity=0.5 ] (327.37,88.4) -- (333.6,87.88) -- (333.6,95.69) -- (323.55,95.69) -- cycle ;
\draw    (276.76,109.86) .. controls (307.97,111.24) and (306.26,143.01) .. (290.69,149.41) ;
\draw  [fill={rgb, 255:red, 0; green, 0; blue, 0 }  ,fill opacity=1 ] (247.23,150.19) -- (269.86,113.75) -- (290.69,149.41) -- cycle ;
\draw    (269.86,113.75) -- (265.7,106.61) ;
\draw    (276.76,109.86) -- (272.6,102.73) ;
\draw  [draw opacity=0][fill={rgb, 255:red, 74; green, 221; blue, 226 }  ,fill opacity=0.5 ] (269.25,111.66) -- (265.7,106.61) -- (272.6,102.73) -- (277.59,111.28) -- cycle ;
\draw    (313.81,145.91) .. controls (298.67,119.07) and (327.13,103.88) .. (340.84,113.55) ;
\draw  [fill={rgb, 255:red, 0; green, 0; blue, 0 }  ,fill opacity=1 ] (364.26,149.51) -- (320.79,149.65) -- (340.84,113.55) -- cycle ;
\draw    (320.79,149.65) -- (316.79,156.87) ;
\draw    (313.81,145.91) -- (309.8,153.13) ;
\draw  [draw opacity=0][fill={rgb, 255:red, 74; green, 221; blue, 226 }  ,fill opacity=0.5 ] (319.3,151.25) -- (316.79,156.87) -- (309.8,153.13) -- (314.61,144.47) -- cycle ;
\draw    (265.7,106.61) -- (272.6,102.73) ;
\draw    (309.8,153.13) -- (316.83,156.79) ;
\draw    (333.6,87.88) -- (333.6,95.69) ;
\draw [color={rgb, 255:red, 244; green, 5; blue, 5 }  ,draw opacity=1 ]   (364.25,102) -- (307.4,74.53) ;
\draw [shift={(304.7,73.22)}, rotate = 25.79] [fill={rgb, 255:red, 244; green, 5; blue, 5 }  ,fill opacity=1 ][line width=0.08]  [draw opacity=0] (3.57,-1.72) -- (0,0) -- (3.57,1.72) -- cycle    ;
\draw [color={rgb, 255:red, 244; green, 5; blue, 5 }  ,draw opacity=1 ]   (364.25,102) -- (342.14,135.74) ;
\draw [shift={(340.5,138.25)}, rotate = 303.23] [fill={rgb, 255:red, 244; green, 5; blue, 5 }  ,fill opacity=1 ][line width=0.08]  [draw opacity=0] (3.57,-1.72) -- (0,0) -- (3.57,1.72) -- cycle    ;
\draw [color={rgb, 255:red, 244; green, 5; blue, 5 }  ,draw opacity=1 ][line width=0.75]    (364.25,102) -- (278.59,137.71) ;
\draw [shift={(275.82,138.86)}, rotate = 337.37] [fill={rgb, 255:red, 244; green, 5; blue, 5 }  ,fill opacity=1 ][line width=0.08]  [draw opacity=0] (3.57,-1.72) -- (0,0) -- (3.57,1.72) -- cycle    ;
\draw  [fill={rgb, 255:red, 155; green, 155; blue, 155 }  ,fill opacity=1 ] (291.73,129.28) .. controls (290.01,130.29) and (287.8,129.69) .. (286.81,127.94) .. controls (285.81,126.19) and (286.4,123.94) .. (288.13,122.93) .. controls (289.85,121.92) and (292.06,122.52) .. (293.05,124.27) .. controls (294.05,126.03) and (293.46,128.27) .. (291.73,129.28) -- cycle ;
\draw  [fill={rgb, 255:red, 155; green, 155; blue, 155 }  ,fill opacity=1 ] (323.42,123.13) .. controls (325.15,124.13) and (325.76,126.37) .. (324.78,128.13) .. controls (323.8,129.89) and (321.6,130.51) .. (319.87,129.52) .. controls (318.14,128.52) and (317.52,126.28) .. (318.5,124.52) .. controls (319.48,122.76) and (321.68,122.14) .. (323.42,123.13) -- cycle ;
\draw [fill={rgb, 255:red, 155; green, 155; blue, 155 }  ,fill opacity=1 ]   (311.3,117.5) -- (319.88,123.4) ;
\draw [fill={rgb, 255:red, 155; green, 155; blue, 155 }  ,fill opacity=1 ]   (308.66,121.92) -- (318.11,126.59) ;
\draw  [fill={rgb, 255:red, 155; green, 155; blue, 155 }  ,fill opacity=1 ] (312.03,120.6) -- (299.83,120.86) -- (305.96,109.84) -- cycle ;
\draw [fill={rgb, 255:red, 155; green, 155; blue, 155 }  ,fill opacity=1 ]   (319.88,123.4) -- (321.64,126.32) ;
\draw [fill={rgb, 255:red, 155; green, 155; blue, 155 }  ,fill opacity=1 ]   (321.64,126.32) -- (318.11,126.59) ;
\draw [fill={rgb, 255:red, 155; green, 155; blue, 155 }  ,fill opacity=1 ]   (311.3,117.5) -- (308.61,112.51) ;
\draw [fill={rgb, 255:red, 155; green, 155; blue, 155 }  ,fill opacity=1 ]   (308.66,121.92) -- (302.94,121.7) ;
\draw [fill={rgb, 255:red, 155; green, 155; blue, 155 }  ,fill opacity=1 ]   (321.64,126.32) -- (312.02,120.6) ;

\draw  [draw opacity=0][fill={rgb, 255:red, 155; green, 155; blue, 155 }  ,fill opacity=1 ] (321.64,126.32) -- (318.11,126.59) -- (308.66,121.92) -- (302.66,121.52) -- (293.3,126.05) -- (289.93,126.11) -- (299.83,120.86) -- (311.86,120.73) -- cycle ;
\draw  [fill={rgb, 255:red, 155; green, 155; blue, 155 }  ,fill opacity=1 ] (302.24,98.67) .. controls (302.24,96.64) and (303.86,95) .. (305.85,95) .. controls (307.84,95) and (309.45,96.64) .. (309.45,98.67) .. controls (309.45,100.69) and (307.84,102.33) .. (305.85,102.33) .. controls (303.86,102.33) and (302.24,100.69) .. (302.24,98.67) -- cycle ;
\draw [fill={rgb, 255:red, 155; green, 155; blue, 155 }  ,fill opacity=1 ]   (303.39,112.17) -- (304.21,101.67) ;
\draw [fill={rgb, 255:red, 155; green, 155; blue, 155 }  ,fill opacity=1 ]   (308.47,112.33) -- (307.82,101.67) ;
\draw  [fill={rgb, 255:red, 155; green, 155; blue, 155 }  ,fill opacity=1 ] (305.68,110) -- (311.91,120.67) -- (299.46,120.67) -- cycle ;
\draw [fill={rgb, 255:red, 155; green, 155; blue, 155 }  ,fill opacity=1 ]   (304.21,101.67) -- (305.85,98.67) ;
\draw [fill={rgb, 255:red, 155; green, 155; blue, 155 }  ,fill opacity=1 ]   (305.85,98.67) -- (307.82,101.67) ;
\draw [fill={rgb, 255:red, 155; green, 155; blue, 155 }  ,fill opacity=1 ]   (303.39,112.17) -- (300.44,117) ;
\draw [fill={rgb, 255:red, 155; green, 155; blue, 155 }  ,fill opacity=1 ]   (308.47,112.33) -- (311.09,117.5) ;
\draw [fill={rgb, 255:red, 155; green, 155; blue, 155 }  ,fill opacity=1 ]   (305.85,98.67) -- (305.68,110) ;
\draw [fill={rgb, 255:red, 155; green, 155; blue, 155 }  ,fill opacity=1 ]   (302.66,121.52) -- (293.3,126.05) ;
\draw [fill={rgb, 255:red, 155; green, 155; blue, 155 }  ,fill opacity=1 ]   (300.26,116.96) -- (291.5,122.87) ;
\draw  [fill={rgb, 255:red, 155; green, 155; blue, 155 }  ,fill opacity=1 ] (299.99,120.65) -- (305.96,109.83) -- (312.18,120.8) -- cycle ;
\draw [fill={rgb, 255:red, 155; green, 155; blue, 155 }  ,fill opacity=1 ]   (293.3,126.05) -- (289.93,126.11) ;
\draw [fill={rgb, 255:red, 155; green, 155; blue, 155 }  ,fill opacity=1 ]   (289.93,126.11) -- (291.5,122.87) ;
\draw [fill={rgb, 255:red, 155; green, 155; blue, 155 }  ,fill opacity=1 ]   (302.66,121.52) -- (308.25,121.7) ;
\draw [fill={rgb, 255:red, 155; green, 155; blue, 155 }  ,fill opacity=1 ]   (300.26,116.96) -- (303.35,112.07) ;
\draw [fill={rgb, 255:red, 155; green, 155; blue, 155 }  ,fill opacity=1 ]   (289.93,126.11) -- (299.66,120.58) ;
\draw  [fill={rgb, 255:red, 155; green, 155; blue, 155 }  ,fill opacity=1 ] (305.24,116.93) .. controls (305.24,116.64) and (305.49,116.4) .. (305.8,116.4) .. controls (306.1,116.4) and (306.35,116.64) .. (306.35,116.93) .. controls (306.35,117.22) and (306.1,117.46) .. (305.8,117.46) .. controls (305.49,117.46) and (305.24,117.22) .. (305.24,116.93) -- cycle ;
\draw  [draw opacity=0][fill={rgb, 255:red, 155; green, 155; blue, 155 }  ,fill opacity=1 ] (305.64,98.75) -- (305.75,109.92) -- (299.46,120.67) -- (289.72,126.19) -- (291.3,122.96) -- (300.05,117.05) -- (303.14,112.15) -- (304,101.75) -- cycle ;
\draw  [draw opacity=0][fill={rgb, 255:red, 155; green, 155; blue, 155 }  ,fill opacity=1 ] (305.85,98.67) -- (307.82,101.67) -- (308.61,112.51) -- (311.3,117.5) -- (319.88,123.4) -- (321.64,126.32) -- (311.86,120.73) -- (305.96,109.83) -- cycle ;
\draw    (311.91,120.67) -- (321.64,126.32) ;
\draw    (305.63,109.77) -- (305.85,98.67) ;

\draw [color={rgb, 255:red, 244; green, 5; blue, 5 }  ,draw opacity=1 ]   (249,87) -- (302.59,115.52) ;
\draw [shift={(305.24,116.93)}, rotate = 208.02] [fill={rgb, 255:red, 244; green, 5; blue, 5 }  ,fill opacity=1 ][line width=0.08]  [draw opacity=0] (3.57,-1.72) -- (0,0) -- (3.57,1.72) -- cycle    ;
\draw  [color={rgb, 255:red, 23; green, 28; blue, 242 }  ,draw opacity=1 ] (307.25,72.75) -- (274.75,72.75) -- (274.75,40.25) ;
\draw   (333.75,40.25) -- (333.75,72.75) -- (301.25,72.75) ;

\draw  [color={rgb, 255:red, 23; green, 28; blue, 242 }  ,draw opacity=1 ] (263.3,135.95) -- (280.86,163.3) -- (253.51,180.85) ;
\draw   (221.64,131.2) -- (248.99,113.65) -- (266.55,141) ;

\draw  [color={rgb, 255:red, 23; green, 28; blue, 242 }  ,draw opacity=1 ] (343.68,140.97) -- (359.93,112.83) -- (388.07,129.08) ;
\draw   (358.57,180.17) -- (330.43,163.92) -- (346.68,135.78) ;

\draw   (461,115.85) .. controls (461,114.69) and (461.94,113.75) .. (463.1,113.75) -- (473.15,113.75) .. controls (474.31,113.75) and (475.25,114.69) .. (475.25,115.85) -- (475.25,122.15) .. controls (475.25,123.31) and (474.31,124.25) .. (473.15,124.25) -- (463.1,124.25) .. controls (461.94,124.25) and (461,123.31) .. (461,122.15) -- cycle ;
\draw    (468,114) -- (468,106) ;
\draw    (467.75,135.5) -- (467.75,124.5) ;
\draw  [draw opacity=0][fill={rgb, 255:red, 155; green, 155; blue, 155 }  ,fill opacity=0.7 ] (463.9,135) -- (473,135) -- (469.1,139.25) -- (460,139.25) -- cycle ;
\draw   (464.5,102.5) .. controls (464.5,100.57) and (466.07,99) .. (468,99) .. controls (469.93,99) and (471.5,100.57) .. (471.5,102.5) .. controls (471.5,104.43) and (469.93,106) .. (468,106) .. controls (466.07,106) and (464.5,104.43) .. (464.5,102.5) -- cycle ;
\draw   (496,115.6) .. controls (496,114.44) and (496.94,113.5) .. (498.1,113.5) -- (508.15,113.5) .. controls (509.31,113.5) and (510.25,114.44) .. (510.25,115.6) -- (510.25,121.9) .. controls (510.25,123.06) and (509.31,124) .. (508.15,124) -- (498.1,124) .. controls (496.94,124) and (496,123.06) .. (496,121.9) -- cycle ;
\draw    (503,113.75) -- (503,105.75) ;
\draw    (502.75,135.25) -- (502.75,124.25) ;
\draw  [draw opacity=0][fill={rgb, 255:red, 155; green, 155; blue, 155 }  ,fill opacity=0.7 ] (498.9,134.75) -- (508,134.75) -- (504.1,139) -- (495,139) -- cycle ;
\draw   (499.5,102.25) .. controls (499.5,100.32) and (501.07,98.75) .. (503,98.75) .. controls (504.93,98.75) and (506.5,100.32) .. (506.5,102.25) .. controls (506.5,104.18) and (504.93,105.75) .. (503,105.75) .. controls (501.07,105.75) and (499.5,104.18) .. (499.5,102.25) -- cycle ;

\draw   (530.25,115.35) .. controls (530.25,114.19) and (531.19,113.25) .. (532.35,113.25) -- (542.4,113.25) .. controls (543.56,113.25) and (544.5,114.19) .. (544.5,115.35) -- (544.5,121.65) .. controls (544.5,122.81) and (543.56,123.75) .. (542.4,123.75) -- (532.35,123.75) .. controls (531.19,123.75) and (530.25,122.81) .. (530.25,121.65) -- cycle ;
\draw    (537.25,113.5) -- (537.25,105.5) ;
\draw    (537,135) -- (537,124) ;
\draw  [draw opacity=0][fill={rgb, 255:red, 155; green, 155; blue, 155 }  ,fill opacity=0.7 ] (533.15,134.5) -- (542.25,134.5) -- (538.35,138.75) -- (529.25,138.75) -- cycle ;
\draw   (533.75,102) .. controls (533.75,100.07) and (535.32,98.5) .. (537.25,98.5) .. controls (539.18,98.5) and (540.75,100.07) .. (540.75,102) .. controls (540.75,103.93) and (539.18,105.5) .. (537.25,105.5) .. controls (535.32,105.5) and (533.75,103.93) .. (533.75,102) -- cycle ;

\draw    (468,99) -- (537.25,98.5) ;
\draw [color={rgb, 255:red, 244; green, 5; blue, 5 }  ,draw opacity=1 ]   (467.48,97.84) -- (472,79.75) ;
\draw [shift={(466.75,100.75)}, rotate = 284.04] [fill={rgb, 255:red, 244; green, 5; blue, 5 }  ,fill opacity=1 ][line width=0.08]  [draw opacity=0] (3.57,-1.72) -- (0,0) -- (3.57,1.72) -- cycle    ;
\draw [color={rgb, 255:red, 244; green, 5; blue, 5 }  ,draw opacity=1 ]   (497,95.15) -- (503,84.75) ;
\draw [shift={(495.5,97.75)}, rotate = 299.98] [fill={rgb, 255:red, 244; green, 5; blue, 5 }  ,fill opacity=1 ][line width=0.08]  [draw opacity=0] (3.57,-1.72) -- (0,0) -- (3.57,1.72) -- cycle    ;
\draw [color={rgb, 255:red, 244; green, 5; blue, 5 }  ,draw opacity=1 ]   (535.41,120.76) -- (510,150) ;
\draw [shift={(537.38,118.5)}, rotate = 130.99] [fill={rgb, 255:red, 244; green, 5; blue, 5 }  ,fill opacity=1 ][line width=0.08]  [draw opacity=0] (3.57,-1.72) -- (0,0) -- (3.57,1.72) -- cycle    ;
\draw  [fill={rgb, 255:red, 74; green, 221; blue, 226 }  ,fill opacity=0.5 ] (185.25,344.25) .. controls (185.25,344.25) and (185.25,344.25) .. (185.25,344.25) .. controls (185.25,344.25) and (185.25,344.25) .. (185.25,344.25) .. controls (185.25,340.04) and (189.95,336.63) .. (195.75,336.63) .. controls (201.46,336.63) and (206.11,339.94) .. (206.25,344.07) -- cycle ;
\draw  [fill={rgb, 255:red, 74; green, 221; blue, 226 }  ,fill opacity=0.5 ] (145,344) .. controls (145,344) and (145,344) .. (145,344) .. controls (145,339.65) and (153.28,336.13) .. (163.5,336.13) .. controls (173.46,336.13) and (181.59,339.48) .. (181.98,343.68) -- cycle ;
\draw  [fill={rgb, 255:red, 74; green, 221; blue, 226 }  ,fill opacity=0.5 ] (209.25,344) .. controls (209.25,344) and (209.25,344) .. (209.25,344) .. controls (209.25,339.65) and (217.53,336.13) .. (227.75,336.13) .. controls (237.71,336.13) and (245.84,339.48) .. (246.23,343.68) -- cycle ;
\draw [line width=1.5]    (165.75,336.13) -- (175.5,336.38) ;
\draw [line width=1.5]    (191,335.63) -- (200.75,335.88) ;
\draw [line width=1.5]    (215,335.88) -- (224.75,336.13) ;
\draw  [fill={rgb, 255:red, 155; green, 155; blue, 155 }  ,fill opacity=1 ] (186.25,319.63) -- (205,319.63) -- (220,334.88) -- (215,334.63) -- (199.5,322.13) -- (191.5,322.13) -- (176,334.88) -- (171,335) -- cycle ;
\draw  [fill={rgb, 255:red, 155; green, 155; blue, 155 }  ,fill opacity=1 ] (191.5,322.13) -- (199.5,322.13) -- (198.5,334.38) -- (192.5,334.38) -- cycle ;
\draw    (193.75,321.88) -- (194.5,334.13) ;
\draw    (197,321.88) -- (196.25,334.38) ;
\draw  [fill={rgb, 255:red, 155; green, 155; blue, 155 }  ,fill opacity=1 ] (194.25,218.88) -- (197,218.88) -- (196.75,319.38) -- (194,319.38) -- cycle ;

\draw [color={rgb, 255:red, 23; green, 28; blue, 242 }  ,draw opacity=1 ]   (145,344) .. controls (140.75,344.63) and (148.75,335.13) .. (161.25,335.38) ;
\draw [color={rgb, 255:red, 23; green, 28; blue, 242 }  ,draw opacity=1 ]   (145,345) -- (160,344.88) ;
\draw [color={rgb, 255:red, 23; green, 28; blue, 242 }  ,draw opacity=1 ]   (246.23,343.68) .. controls (244.25,345.63) and (250,338.88) .. (234,335.38) ;
\draw [color={rgb, 255:red, 23; green, 28; blue, 242 }  ,draw opacity=1 ]   (231.23,344.8) -- (246.23,344.68) ;
\draw [line width=1.5]    (360.75,324.46) -- (370.5,324.71) ;
\draw [line width=1.5]    (386,323.96) -- (395.75,324.21) ;
\draw [line width=1.5]    (410,324.21) -- (419.75,324.46) ;
\draw  [fill={rgb, 255:red, 155; green, 155; blue, 155 }  ,fill opacity=1 ] (381.25,307.96) -- (400,307.96) -- (415,323.21) -- (410,322.96) -- (394.5,310.46) -- (386.5,310.46) -- (371,323.21) -- (366,323.33) -- cycle ;
\draw  [fill={rgb, 255:red, 155; green, 155; blue, 155 }  ,fill opacity=1 ] (386.5,310.46) -- (394.5,310.46) -- (393.5,322.71) -- (387.5,322.71) -- cycle ;
\draw    (388.75,310.21) -- (389.5,322.46) ;
\draw    (392,310.21) -- (391.25,322.71) ;
\draw  [fill={rgb, 255:red, 155; green, 155; blue, 155 }  ,fill opacity=1 ] (389.25,207.21) -- (392,207.21) -- (391.75,307.71) -- (389,307.71) -- cycle ;

\draw  [fill={rgb, 255:red, 74; green, 221; blue, 226 }  ,fill opacity=0.5 ] (330.67,344.17) .. controls (340.33,344.17) and (354.67,323.17) .. (367,326.83) .. controls (379.33,330.5) and (378,342.83) .. (377,343.83) .. controls (376,344.83) and (375.33,344.5) .. (360.33,344.5) .. controls (345.33,344.5) and (321,344.17) .. (330.67,344.17) -- cycle ;
\draw  [fill={rgb, 255:red, 74; green, 221; blue, 226 }  ,fill opacity=0.5 ] (450.67,343.5) .. controls (438.33,342.5) and (426,324.5) .. (415.67,325.17) .. controls (405.33,325.83) and (400.33,344.17) .. (405,343.83) .. controls (409.67,343.5) and (418,344.17) .. (424.67,344.17) .. controls (431.33,344.17) and (463,344.5) .. (450.67,343.5) -- cycle ;
\draw  [fill={rgb, 255:red, 74; green, 221; blue, 226 }  ,fill opacity=0.5 ] (381.48,343.22) .. controls (380.9,341.63) and (380.58,339.87) .. (380.58,338.02) .. controls (380.58,330.92) and (385.28,325.17) .. (391.08,325.17) .. controls (396.88,325.17) and (401.58,330.92) .. (401.58,338.02) .. controls (401.58,339.85) and (401.27,341.59) .. (400.71,343.17) -- cycle ;
\draw  [fill={rgb, 255:red, 244; green, 5; blue, 5 }  ,fill opacity=0.5 ] (362,334.37) -- (365.67,327.83) -- (369.33,334.37) -- (367.5,334.37) -- (367.5,344.17) -- (363.83,344.17) -- (363.83,334.37) -- cycle ;
\draw  [fill={rgb, 255:red, 244; green, 5; blue, 5 }  ,fill opacity=0.5 ] (387.67,332.7) -- (391.33,326.17) -- (395,332.7) -- (393.17,332.7) -- (393.17,342.5) -- (389.5,342.5) -- (389.5,332.7) -- cycle ;
\draw  [fill={rgb, 255:red, 244; green, 5; blue, 5 }  ,fill opacity=0.5 ] (412,333.7) -- (415.67,327.17) -- (419.33,333.7) -- (417.5,333.7) -- (417.5,343.5) -- (413.83,343.5) -- (413.83,333.7) -- cycle ;
\draw [color={rgb, 255:red, 155; green, 155; blue, 155 }  ,draw opacity=1 ]   (220,334.88) -- (307,334.63) ;
\draw [color={rgb, 255:red, 155; green, 155; blue, 155 }  ,draw opacity=1 ]   (274,327.08) -- (361,326.83) ;
\draw    (303.67,329.5) -- (303.67,331.83) ;
\draw [shift={(303.67,334.83)}, rotate = 270] [fill={rgb, 255:red, 0; green, 0; blue, 0 }  ][line width=0.08]  [draw opacity=0] (3.57,-1.72) -- (0,0) -- (3.57,1.72) -- cycle    ;
\draw [shift={(303.67,326.5)}, rotate = 90] [fill={rgb, 255:red, 0; green, 0; blue, 0 }  ][line width=0.08]  [draw opacity=0] (3.57,-1.72) -- (0,0) -- (3.57,1.72) -- cycle    ;
\draw [color={rgb, 255:red, 244; green, 5; blue, 5 }  ,draw opacity=1 ]   (196.76,217.24) -- (209.67,208.83) ;
\draw [shift={(194.25,218.88)}, rotate = 326.92] [fill={rgb, 255:red, 244; green, 5; blue, 5 }  ,fill opacity=1 ][line width=0.08]  [draw opacity=0] (3.57,-1.72) -- (0,0) -- (3.57,1.72) -- cycle    ;
\draw [color={rgb, 255:red, 244; green, 5; blue, 5 }  ,draw opacity=1 ]   (197.76,259.9) -- (210.67,251.5) ;
\draw [shift={(195.25,261.54)}, rotate = 326.92] [fill={rgb, 255:red, 244; green, 5; blue, 5 }  ,fill opacity=1 ][line width=0.08]  [draw opacity=0] (3.57,-1.72) -- (0,0) -- (3.57,1.72) -- cycle    ;
\draw [color={rgb, 255:red, 155; green, 155; blue, 155 }  ,draw opacity=1 ]   (389.25,207.21) -- (429,207.21) ;
\draw [color={rgb, 255:red, 155; green, 155; blue, 155 }  ,draw opacity=1 ]   (391.25,322.71) -- (431,322.71) ;
\draw [color={rgb, 255:red, 155; green, 155; blue, 155 }  ,draw opacity=1 ]   (409.25,251.75) -- (409.13,210.21) ;
\draw [shift={(409.13,207.21)}, rotate = 89.84] [fill={rgb, 255:red, 155; green, 155; blue, 155 }  ,fill opacity=1 ][line width=0.08]  [draw opacity=0] (8.93,-4.29) -- (0,0) -- (8.93,4.29) -- cycle    ;
\draw [color={rgb, 255:red, 155; green, 155; blue, 155 }  ,draw opacity=1 ]   (409.37,319.29) -- (409.25,278.75) ;
\draw [shift={(409.38,322.29)}, rotate = 269.84] [fill={rgb, 255:red, 155; green, 155; blue, 155 }  ,fill opacity=1 ][line width=0.08]  [draw opacity=0] (8.93,-4.29) -- (0,0) -- (8.93,4.29) -- cycle    ;
\draw [color={rgb, 255:red, 23; green, 28; blue, 242 }  ,draw opacity=1 ][line width=1.5]    (331.75,343.25) .. controls (336.25,341.75) and (341.25,338.75) .. (348,334.63) ;
\draw [color={rgb, 255:red, 23; green, 28; blue, 242 }  ,draw opacity=1 ][line width=1.5]    (330.67,344.17) -- (345.67,344.04) ;
\draw [color={rgb, 255:red, 23; green, 28; blue, 242 }  ,draw opacity=1 ][line width=1.5]    (450.67,343.5) .. controls (449.25,342.75) and (443.25,342.25) .. (434.25,333.75) ;
\draw [color={rgb, 255:red, 23; green, 28; blue, 242 }  ,draw opacity=1 ][line width=1.5]    (435.67,343.63) -- (450.67,343.5) ;
\draw  [color={rgb, 255:red, 244; green, 5; blue, 5 }  ,draw opacity=1 ] (191.25,335.88) .. controls (191.25,333.25) and (193.38,331.13) .. (196,331.13) .. controls (198.62,331.13) and (200.75,333.25) .. (200.75,335.88) .. controls (200.75,338.5) and (198.62,340.63) .. (196,340.63) .. controls (193.38,340.63) and (191.25,338.5) .. (191.25,335.88) -- cycle ; \draw  [color={rgb, 255:red, 244; green, 5; blue, 5 }  ,draw opacity=1 ] (192.64,332.52) -- (199.36,339.23) ; \draw  [color={rgb, 255:red, 244; green, 5; blue, 5 }  ,draw opacity=1 ] (199.36,332.52) -- (192.64,339.23) ;
\draw  [color={rgb, 255:red, 244; green, 5; blue, 5 }  ,draw opacity=1 ] (386.58,334.75) .. controls (386.58,332.13) and (388.71,330) .. (391.33,330) .. controls (393.96,330) and (396.08,332.13) .. (396.08,334.75) .. controls (396.08,337.37) and (393.96,339.5) .. (391.33,339.5) .. controls (388.71,339.5) and (386.58,337.37) .. (386.58,334.75) -- cycle ; \draw  [color={rgb, 255:red, 244; green, 5; blue, 5 }  ,draw opacity=1 ] (387.97,331.39) -- (394.69,338.11) ; \draw  [color={rgb, 255:red, 244; green, 5; blue, 5 }  ,draw opacity=1 ] (394.69,331.39) -- (387.97,338.11) ;
\draw [color={rgb, 255:red, 248; green, 231; blue, 28 }  ,draw opacity=1 ]   (196,335.88) -- (195.51,268.37) ;
\draw [shift={(195.5,266.38)}, rotate = 89.59] [color={rgb, 255:red, 248; green, 231; blue, 28 }  ,draw opacity=1 ][line width=0.75]    (10.93,-3.29) .. controls (6.95,-1.4) and (3.31,-0.3) .. (0,0) .. controls (3.31,0.3) and (6.95,1.4) .. (10.93,3.29)   ;
\draw [color={rgb, 255:red, 91; green, 163; blue, 3 }  ,draw opacity=1 ]   (196,335.88) -- (237.5,335.76) ;
\draw [shift={(239.5,335.75)}, rotate = 179.84] [color={rgb, 255:red, 91; green, 163; blue, 3 }  ,draw opacity=1 ][line width=0.75]    (10.93,-3.29) .. controls (6.95,-1.4) and (3.31,-0.3) .. (0,0) .. controls (3.31,0.3) and (6.95,1.4) .. (10.93,3.29)   ;
\draw [color={rgb, 255:red, 244; green, 5; blue, 5 }  ,draw opacity=1 ]   (197.52,338.46) -- (222,380) ;
\draw [shift={(196,335.88)}, rotate = 59.49] [fill={rgb, 255:red, 244; green, 5; blue, 5 }  ,fill opacity=1 ][line width=0.08]  [draw opacity=0] (3.57,-1.72) -- (0,0) -- (3.57,1.72) -- cycle    ;
\draw [color={rgb, 255:red, 248; green, 231; blue, 28 }  ,draw opacity=1 ]   (391.5,334.88) -- (391.01,267.37) ;
\draw [shift={(391,265.38)}, rotate = 89.59] [color={rgb, 255:red, 248; green, 231; blue, 28 }  ,draw opacity=1 ][line width=0.75]    (10.93,-3.29) .. controls (6.95,-1.4) and (3.31,-0.3) .. (0,0) .. controls (3.31,0.3) and (6.95,1.4) .. (10.93,3.29)   ;
\draw [color={rgb, 255:red, 91; green, 163; blue, 3 }  ,draw opacity=1 ]   (391.5,334.88) -- (433,334.76) ;
\draw [shift={(435,334.75)}, rotate = 179.84] [color={rgb, 255:red, 91; green, 163; blue, 3 }  ,draw opacity=1 ][line width=0.75]    (10.93,-3.29) .. controls (6.95,-1.4) and (3.31,-0.3) .. (0,0) .. controls (3.31,0.3) and (6.95,1.4) .. (10.93,3.29)   ;
\draw [color={rgb, 255:red, 244; green, 5; blue, 5 }  ,draw opacity=1 ]   (160.65,335.44) -- (155,308) ;
\draw [shift={(161.25,338.38)}, rotate = 258.37] [fill={rgb, 255:red, 244; green, 5; blue, 5 }  ,fill opacity=1 ][line width=0.08]  [draw opacity=0] (3.57,-1.72) -- (0,0) -- (3.57,1.72) -- cycle    ;
\draw [color={rgb, 255:red, 244; green, 5; blue, 5 }  ,draw opacity=1 ]   (366.48,332.96) -- (440.5,293.75) ;
\draw [shift={(363.83,334.37)}, rotate = 332.09] [fill={rgb, 255:red, 244; green, 5; blue, 5 }  ,fill opacity=1 ][line width=0.08]  [draw opacity=0] (3.57,-1.72) -- (0,0) -- (3.57,1.72) -- cycle    ;
\draw [color={rgb, 255:red, 244; green, 5; blue, 5 }  ,draw opacity=1 ]   (393.64,332.83) -- (440.5,293.75) ;
\draw [shift={(391.33,334.75)}, rotate = 320.18] [fill={rgb, 255:red, 244; green, 5; blue, 5 }  ,fill opacity=1 ][line width=0.08]  [draw opacity=0] (3.57,-1.72) -- (0,0) -- (3.57,1.72) -- cycle    ;
\draw [color={rgb, 255:red, 244; green, 5; blue, 5 }  ,draw opacity=1 ]   (417.49,334.14) -- (440.5,293.75) ;
\draw [shift={(416,336.75)}, rotate = 299.67] [fill={rgb, 255:red, 244; green, 5; blue, 5 }  ,fill opacity=1 ][line width=0.08]  [draw opacity=0] (3.57,-1.72) -- (0,0) -- (3.57,1.72) -- cycle    ;
\draw [color={rgb, 255:red, 244; green, 5; blue, 5 }  ,draw opacity=1 ]   (344.06,339.55) -- (359.5,380.25) ;
\draw [shift={(343,336.75)}, rotate = 69.23] [fill={rgb, 255:red, 244; green, 5; blue, 5 }  ,fill opacity=1 ][line width=0.08]  [draw opacity=0] (3.57,-1.72) -- (0,0) -- (3.57,1.72) -- cycle    ;
\draw [color={rgb, 255:red, 244; green, 5; blue, 5 }  ,draw opacity=1 ]   (336.72,347.71) -- (359.5,380.25) ;
\draw [shift={(335,345.25)}, rotate = 55.01] [fill={rgb, 255:red, 244; green, 5; blue, 5 }  ,fill opacity=1 ][line width=0.08]  [draw opacity=0] (3.57,-1.72) -- (0,0) -- (3.57,1.72) -- cycle    ;
\draw [color={rgb, 255:red, 244; green, 5; blue, 5 }  ,draw opacity=1 ]   (172.64,336.32) -- (232,304.25) ;
\draw [shift={(170,337.75)}, rotate = 331.62] [fill={rgb, 255:red, 244; green, 5; blue, 5 }  ,fill opacity=1 ][line width=0.08]  [draw opacity=0] (3.57,-1.72) -- (0,0) -- (3.57,1.72) -- cycle    ;
\draw [color={rgb, 255:red, 244; green, 5; blue, 5 }  ,draw opacity=1 ]   (198.25,333.9) -- (232,304.25) ;
\draw [shift={(196,335.88)}, rotate = 318.7] [fill={rgb, 255:red, 244; green, 5; blue, 5 }  ,fill opacity=1 ][line width=0.08]  [draw opacity=0] (3.57,-1.72) -- (0,0) -- (3.57,1.72) -- cycle    ;
\draw [color={rgb, 255:red, 244; green, 5; blue, 5 }  ,draw opacity=1 ]   (225.42,333.2) -- (232,304.25) ;
\draw [shift={(224.75,336.13)}, rotate = 282.81] [fill={rgb, 255:red, 244; green, 5; blue, 5 }  ,fill opacity=1 ][line width=0.08]  [draw opacity=0] (3.57,-1.72) -- (0,0) -- (3.57,1.72) -- cycle    ;
\draw  [color={rgb, 255:red, 74; green, 221; blue, 226 }  ,draw opacity=1 ][fill={rgb, 255:red, 74; green, 221; blue, 226 }  ,fill opacity=0.2 ] (213.33,479.69) -- (303,479.69) -- (303,493.02) -- (213.33,493.02) -- cycle ;
\draw  [color={rgb, 255:red, 0; green, 0; blue, 0 }  ,draw opacity=1 ][fill={rgb, 255:red, 74; green, 221; blue, 226 }  ,fill opacity=0.2 ] (388.33,446.02) -- (490.33,450.02) -- (403.33,493.02) -- (304,493.02) -- cycle ;
\draw    (188.33,493.02) -- (188.33,427.69) ;
\draw [shift={(188.33,425.69)}, rotate = 90] [color={rgb, 255:red, 0; green, 0; blue, 0 }  ][line width=0.75]    (10.93,-3.29) .. controls (6.95,-1.4) and (3.31,-0.3) .. (0,0) .. controls (3.31,0.3) and (6.95,1.4) .. (10.93,3.29)   ;
\draw    (188.33,493.02) -- (531.33,493.68) ;
\draw [shift={(533.33,493.69)}, rotate = 180.11] [color={rgb, 255:red, 0; green, 0; blue, 0 }  ][line width=0.75]    (10.93,-3.29) .. controls (6.95,-1.4) and (3.31,-0.3) .. (0,0) .. controls (3.31,0.3) and (6.95,1.4) .. (10.93,3.29)   ;
\draw  [color={rgb, 255:red, 0; green, 0; blue, 0 }  ,draw opacity=1 ] (188.33,492.02) -- (213.33,492.02) -- (213.33,493.02) -- (188.33,493.02) -- cycle ;
\draw  [draw opacity=0] (360.16,467.25) .. controls (360.42,467.65) and (360.68,468.07) .. (360.94,468.49) .. controls (364.32,474.2) and (365.64,480.62) .. (365.13,487.01) -- (328.1,487.98) -- cycle ; \draw   (360.16,467.25) .. controls (360.42,467.65) and (360.68,468.07) .. (360.94,468.49) .. controls (364.32,474.2) and (365.64,480.62) .. (365.13,487.01) ;  
\draw [line width=2.25]    (213.33,479.69) -- (303,479.69) ;
\draw  [draw opacity=0] (458.72,474.5) .. controls (460.96,479.42) and (461.79,484.75) .. (461.37,490.05) -- (424.33,491.02) -- cycle ; \draw   (458.72,474.5) .. controls (460.96,479.42) and (461.79,484.75) .. (461.37,490.05) ;  
\draw [color={rgb, 255:red, 155; green, 155; blue, 155 }  ,draw opacity=1 ]   (191,504) -- (300.8,504) ;
\draw [shift={(303.8,504)}, rotate = 180] [fill={rgb, 255:red, 155; green, 155; blue, 155 }  ,fill opacity=1 ][line width=0.08]  [draw opacity=0] (8.93,-4.29) -- (0,0) -- (8.93,4.29) -- cycle    ;
\draw [shift={(188,504)}, rotate = 0] [fill={rgb, 255:red, 155; green, 155; blue, 155 }  ,fill opacity=1 ][line width=0.08]  [draw opacity=0] (8.93,-4.29) -- (0,0) -- (8.93,4.29) -- cycle    ;
\draw [color={rgb, 255:red, 155; green, 155; blue, 155 }  ,draw opacity=1 ]   (188.33,437.02) -- (188.8,494) ;
\draw [color={rgb, 255:red, 155; green, 155; blue, 155 }  ,draw opacity=1 ]   (303,494.02) -- (303.47,520) ;
\draw [color={rgb, 255:red, 155; green, 155; blue, 155 }  ,draw opacity=1 ]   (213.33,436.02) -- (213.8,493) ;
\draw [color={rgb, 255:red, 155; green, 155; blue, 155 }  ,draw opacity=1 ]   (190.43,453.51) -- (210.23,453.51) ;
\draw [shift={(213.23,453.51)}, rotate = 180] [fill={rgb, 255:red, 155; green, 155; blue, 155 }  ,fill opacity=1 ][line width=0.08]  [draw opacity=0] (8.93,-4.29) -- (0,0) -- (8.93,4.29) -- cycle    ;
\draw [shift={(187.43,453.51)}, rotate = 0] [fill={rgb, 255:red, 155; green, 155; blue, 155 }  ,fill opacity=1 ][line width=0.08]  [draw opacity=0] (8.93,-4.29) -- (0,0) -- (8.93,4.29) -- cycle    ;
\draw [color={rgb, 255:red, 155; green, 155; blue, 155 }  ,draw opacity=1 ]   (187.77,491.01) -- (188.23,516.99) ;
\draw [line width=1.5]    (502.58,324.32) -- (512.08,326.53) ;
\draw [line width=1.5]    (527.41,328.92) -- (536.91,331.12) ;
\draw [line width=1.5]    (550.87,333.99) -- (560.37,336.2) ;
\draw  [fill={rgb, 255:red, 155; green, 155; blue, 155 }  ,fill opacity=1 ] (525.98,312.29) -- (544.35,316.06) -- (555.97,334.02) -- (551.12,332.77) -- (538.46,317.4) -- (530.62,315.79) -- (512.87,325.16) -- (507.95,324.28) -- cycle ;
\draw  [fill={rgb, 255:red, 155; green, 155; blue, 155 }  ,fill opacity=1 ] (530.62,315.79) -- (538.46,317.4) -- (535.01,329.2) -- (529.14,327.99) -- cycle ;
\draw    (532.88,316) -- (531.14,328.15) ;
\draw    (536.06,316.66) -- (532.81,328.75) ;
\draw  [fill={rgb, 255:red, 155; green, 155; blue, 155 }  ,fill opacity=1 ] (554.09,215.21) -- (556.79,215.76) -- (536.32,314.16) -- (533.62,313.6) -- cycle ;

\draw  [fill={rgb, 255:red, 74; green, 221; blue, 226 }  ,fill opacity=0.5 ] (471.67,342.17) .. controls (481.33,342.17) and (495.67,321.17) .. (508,324.83) .. controls (520.33,328.5) and (519,340.83) .. (518,341.83) .. controls (517,342.83) and (516.33,342.5) .. (501.33,342.5) .. controls (486.33,342.5) and (462,342.17) .. (471.67,342.17) -- cycle ;
\draw  [fill={rgb, 255:red, 244; green, 5; blue, 5 }  ,fill opacity=0.5 ] (503,332.37) -- (506.67,325.83) -- (510.33,332.37) -- (508.5,332.37) -- (508.5,342.17) -- (504.83,342.17) -- (504.83,332.37) -- cycle ;
\draw [color={rgb, 255:red, 23; green, 28; blue, 242 }  ,draw opacity=1 ][line width=1.5]    (469.33,341.83) .. controls (478.33,339.5) and (484,335.83) .. (489,332.63) ;
\draw [color={rgb, 255:red, 23; green, 28; blue, 242 }  ,draw opacity=1 ][line width=1.5]    (469.33,341.83) -- (486.67,342.04) ;
\draw  [fill={rgb, 255:red, 74; green, 221; blue, 226 }  ,fill opacity=0.5 ] (547.25,342.02) .. controls (547.25,342.02) and (547.25,342.02) .. (547.25,342.02) .. controls (547.25,342.02) and (547.25,342.02) .. (547.25,342.02) .. controls (547.25,338.05) and (555.53,334.68) .. (565.75,334.5) .. controls (575.69,334.32) and (583.8,337.22) .. (584.23,341.03) -- cycle ;
\draw [color={rgb, 255:red, 23; green, 28; blue, 242 }  ,draw opacity=1 ]   (584.23,340.68) .. controls (582.25,342.63) and (587.33,338) .. (571.33,334.5) ;
\draw [color={rgb, 255:red, 23; green, 28; blue, 242 }  ,draw opacity=1 ]   (569.23,341.8) -- (584.23,341.68) ;
\draw  [fill={rgb, 255:red, 74; green, 221; blue, 226 }  ,fill opacity=0.5 ] (521.25,341.19) .. controls (521.25,341.19) and (521.25,341.19) .. (521.25,341.19) .. controls (521.25,341.19) and (521.25,341.19) .. (521.25,341.19) .. controls (521.25,335.28) and (525.95,330.5) .. (531.75,330.5) .. controls (537.49,330.5) and (542.15,335.19) .. (542.25,341) -- cycle ;
\draw  [fill={rgb, 255:red, 244; green, 5; blue, 5 }  ,fill opacity=0.5 ] (528.08,335.18) -- (531.75,331.17) -- (535.42,335.18) -- (533.58,335.18) -- (533.58,341.19) -- (529.92,341.19) -- (529.92,335.18) -- cycle ;
\draw  [fill={rgb, 255:red, 244; green, 5; blue, 5 }  ,fill opacity=0.5 ] (554.75,338.04) -- (558.42,336.17) -- (562.08,338.04) -- (560.25,338.04) -- (560.25,340.85) -- (556.58,340.85) -- (556.58,338.04) -- cycle ;
\draw [color={rgb, 255:red, 155; green, 155; blue, 155 }  ,draw opacity=1 ] [dash pattern={on 4.5pt off 4.5pt}]  (529.92,335.18) -- (529.79,201.63) ;
\draw [shift={(529.79,199.63)}, rotate = 89.95] [color={rgb, 255:red, 155; green, 155; blue, 155 }  ,draw opacity=1 ][line width=0.75]    (10.93,-3.29) .. controls (6.95,-1.4) and (3.31,-0.3) .. (0,0) .. controls (3.31,0.3) and (6.95,1.4) .. (10.93,3.29)   ;
\draw [color={rgb, 255:red, 155; green, 155; blue, 155 }  ,draw opacity=1 ] [dash pattern={on 4.5pt off 4.5pt}]  (529.92,335.18) -- (557.59,202.96) ;
\draw [shift={(558,201)}, rotate = 101.82] [color={rgb, 255:red, 155; green, 155; blue, 155 }  ,draw opacity=1 ][line width=0.75]    (10.93,-3.29) .. controls (6.95,-1.4) and (3.31,-0.3) .. (0,0) .. controls (3.31,0.3) and (6.95,1.4) .. (10.93,3.29)   ;

\draw (458,38.75) node [anchor=north west][inner sep=0.75pt]   [align=left] {Spherical Joint\\(S)};
\draw (499.75,67.13) node [anchor=north west][inner sep=0.75pt]   [align=left] {Platform};
\draw (449,146.88) node [anchor=north west][inner sep=0.75pt]   [align=left] {Prismatic Joint (P)};
\draw (292.67,305.83) node [anchor=north west][inner sep=0.75pt]   [align=left] {$\displaystyle h_{i}$};
\draw (205.58,193.88) node [anchor=north west][inner sep=0.75pt]   [align=left] {Tip};
\draw (206.58,236.54) node [anchor=north west][inner sep=0.75pt]   [align=left] {Antenna};
\draw (402.25,254.75) node [anchor=north west][inner sep=0.75pt]   [align=left] {$\displaystyle L$};
\draw (215.58,380.88) node [anchor=north west][inner sep=0.75pt]   [align=left] {$\displaystyle C_{p}$};
\draw (191.08,349.5) node [anchor=north west][inner sep=0.75pt]  [color={rgb, 255:red, 155; green, 155; blue, 155 }  ,opacity=1 ] [align=left] {$\displaystyle X$};
\draw (235.08,348) node [anchor=north west][inner sep=0.75pt]  [color={rgb, 255:red, 155; green, 155; blue, 155 }  ,opacity=1 ] [align=left] {$\displaystyle Y$};
\draw (174.58,253.5) node [anchor=north west][inner sep=0.75pt]  [color={rgb, 255:red, 155; green, 155; blue, 155 }  ,opacity=1 ] [align=left] {$\displaystyle Z$};
\draw (386.58,348.5) node [anchor=north west][inner sep=0.75pt]  [color={rgb, 255:red, 155; green, 155; blue, 155 }  ,opacity=1 ] [align=left] {$\displaystyle X$};
\draw (430.58,347) node [anchor=north west][inner sep=0.75pt]  [color={rgb, 255:red, 155; green, 155; blue, 155 }  ,opacity=1 ] [align=left] {$\displaystyle Y$};
\draw (370.08,252.5) node [anchor=north west][inner sep=0.75pt]  [color={rgb, 255:red, 155; green, 155; blue, 155 }  ,opacity=1 ] [align=left] {$\displaystyle Z$};
\draw (123,268.88) node [anchor=north west][inner sep=0.75pt]   [align=left] {Dielectric\\liquid};
\draw (432,251) node [anchor=north west][inner sep=0.75pt]   [align=left] {Prismatic\\Joint};
\draw (303.5,378) node [anchor=north west][inner sep=0.75pt]   [align=left] {Electrodes};
\draw (203.5,264.5) node [anchor=north west][inner sep=0.75pt]   [align=left] {Passive spherical\\behavior};
\draw (278,175) node [anchor=north west][inner sep=0.75pt]   [align=left] {a)};
\draw (477,175) node [anchor=north west][inner sep=0.75pt]   [align=left] {b)};
\draw (346,406) node [anchor=north west][inner sep=0.75pt]   [align=left] {c)};
\draw (363.75,89.25) node [anchor=north west][inner sep=0.75pt]   [align=left] {Electrodes};
\draw (190.5,42.5) node [anchor=north west][inner sep=0.75pt]   [align=left] {3PS\\parallel \\platform};
\draw (268,23.5) node [anchor=north west][inner sep=0.75pt]  [color={rgb, 255:red, 23; green, 28; blue, 242 }  ,opacity=1 ] [align=left] {+};
\draw (329,22.38) node [anchor=north west][inner sep=0.75pt]   [align=left] {\mbox{-}};
\draw (242,152.63) node [anchor=north west][inner sep=0.75pt]  [color={rgb, 255:red, 23; green, 28; blue, 242 }  ,opacity=1 ] [align=left] {+};
\draw (212.25,118.88) node [anchor=north west][inner sep=0.75pt]   [align=left] {\mbox{-}};
\draw (357,168.63) node [anchor=north west][inner sep=0.75pt]   [align=left] {\mbox{-}};
\draw (385.75,118.63) node [anchor=north west][inner sep=0.75pt]  [color={rgb, 255:red, 23; green, 28; blue, 242 }  ,opacity=1 ] [align=left] {+};
\draw (291,27) node [anchor=north west][inner sep=0.75pt]   [align=left] {$\displaystyle U_{\text{in}}^{1}$};
\draw (374,146) node [anchor=north west][inner sep=0.75pt]   [align=left] {$\displaystyle U_{\text{in}}^{2}$};
\draw (209,147) node [anchor=north west][inner sep=0.75pt]   [align=left] {$\displaystyle U_{\text{in}}^{3}$};
\draw (341,470.02) node [anchor=north west][inner sep=0.75pt]  [font=\normalsize] [align=left] {$\displaystyle \delta _{1}^{i}$};
\draw (437,475.02) node [anchor=north west][inner sep=0.75pt]  [font=\normalsize] [align=left] {$\displaystyle \theta _{i}$};
\draw (479,454.02) node [anchor=north west][inner sep=0.75pt]  [font=\normalsize] [align=left] {$\displaystyle ( 1/2) L_{v}$};
\draw (380,455.72) node [anchor=north west][inner sep=0.75pt]   [align=left] {$\displaystyle A_{s}^{i}$};
\draw (294,433.39) node [anchor=north west][inner sep=0.75pt]  [font=\normalsize] [align=left] {$\displaystyle ( 1/2) l_{p}^{i}$};
\draw (427,424.39) node [anchor=north west][inner sep=0.75pt]  [font=\normalsize] [align=left] {$\displaystyle ( 1/2) l_{p}^{i}$};
\draw (191.57,456.51) node [anchor=north west][inner sep=0.75pt]   [align=left] {$\displaystyle l_{e}^{i}$};
\draw (237.57,503.51) node [anchor=north west][inner sep=0.75pt]   [align=left] {$\displaystyle L_{e}$};
\draw (344,505) node [anchor=north west][inner sep=0.75pt]   [align=left] {d)};

\end{tikzpicture}